\documentclass[review]{elsarticle}

\usepackage{lineno,hyperref}
\usepackage{amsthm}
 \usepackage{amssymb}
 \usepackage{amsmath}
 \usepackage{mathrsfs}
 \usepackage{mathtools}
 \usepackage{enumerate}
 \newtheorem{thm}{Theorem}[section]
 
 \newtheorem{lemma}[thm]{Lemma}
 \newtheorem*{pro}{Proof}
 \newtheorem{remark}{Remark}[section]
 \theoremstyle{definition}
 \newtheorem{definition}[thm]{Definition}
 \theoremstyle{remark}
 
 \numberwithin{equation}{section}
 
\journal{Journal of \LaTeX\ Templates}

\begin{document}

\begin{frontmatter}

\title{Global existence, uniqueness and $L^{\infty}$-bound of weak solutions of fractional time-space Keller-Segel system\tnoteref{mytitlenote}}
\tnotetext[mytitlenote]{This work is supported by the State Key Program of National Natural Science of China under Grant No.91324201. This work is also supported by the Fundamental Research Funds for the Central Universities of China under Grant 2018IB017, Equipment Pre-Research Ministry of Education Joint Fund Grant 6141A02033703 and the Natural Science Foundation of Hubei Province of China under Grant 2014CFB865.}


\author[mymainaddress]{Liujie Guo}
\ead{2252987468@qq.com}

\author[mymainaddress]{Fei Gao\corref{mycorrespondingauthor}}
\cortext[mycorrespondingauthor]{Corresponding author}
\ead{gaof@whut.edu.cn}

\author[mymainaddress]{Hui Zhan}
\ead{2432593867@qq.com}

\address[mymainaddress]{Department of Mathematics and Center for Mathematical Sciences, Wuhan University of Technology, Wuhan, 430070, China}

\begin{abstract}
This paper studies the properties of weak solutions to a class of space-time fractional parabolic-elliptic Keller-Segel equations with logistic source terms in $\mathbb{R}^{n}$, $n\geq 2$. The global existence and $L^{\infty}$-bound of weak solutions are established. We mainly divide the damping coefficient into two cases: (i) $b>1-\frac{\alpha}{n}$, for any initial value and birth rate; (ii) $0<b\leq 1-\frac{\alpha}{n}$, for small initial value and small birth rate. The existence result is obtained by 
verifying the existence of a solution to the constructed regularization equation and incorporate the generalized compactness criterion of time fractional partial differential equation. At the same time, we get the $L^{\infty}$-bound of weak solutions by establishing the fractional differential inequality and using the Moser iterative method. Furthermore, we prove the uniqueness of weak solutions by using the hyper-contractive estimates when the damping coefficient is strong. Finally, we also propose a blow-up criterion for weak solutions, that is, if a weak solution blows up in finite time, then for all $h>q$, the $L^{h}$-norms of the weak solution blow up at the same time.
\end{abstract}


\end{frontmatter}

\section{Introduction}
In this paper, we study the existence, $L^{\infty}$-bound and uniqueness for the following fractional time-space parabolic-elliptic Keller-Segel model:
\begin{equation}\label{1.1}
\left\{
\begin{aligned}
    & \frac{\partial^{\beta}u}{\partial t^{\beta}} = -(-\Delta)^{-\frac{\alpha}{2}}u - \nabla \cdot (u\nabla v) + au-bu^{2},  & x\in \mathbb{R}^{n}, t>0,\\
    & 0=\Delta v+u, & x\in \mathbb{R}^{n}, t>0,\\
    & u(x,0)  = u_0, & x\in \mathbb{R}^{n},\\
\end{aligned}
\right.
\end{equation}
where $\beta\in (0,1)$, $\alpha\in (1,2)$, $n\geq 2$, $a\geq 0$ and $b>0$ represent the birth rate and damping coefficient, respectively. $u_{0}=u_{0}(x): \mathbb{R}^{n}\to \mathbb{R}$ is any given initial value. $\partial_{t}^{\beta}$ is weak Caputo fractional derivative operator of order $\beta$ introduced in \cite{li2018some,li2018generalized}. When function $c$ is absolutely continuous in time, the definition of weak Caputo derivative is reduced to the following traditional form:
\begin{equation}
    \partial_{t}^{\beta}c(t)=\frac{1}{\Gamma(1-\beta)}\int_{0}^{t}(t-s)^{-\beta}\dot{c}(s)ds,
\end{equation}
where $\gamma$ is the Gamma function and $\dot{c}(t)$ is the first order integer derivative of function $c(t)$ with respect to independent variable $t$. The nonlocal operator $(-\Delta)^{\frac{\alpha}{2}}$ known as the Laplacian of order $\frac{\alpha}{2}$, is given by the Fourier multiplier
$$D^{\alpha}u(x):=(-\Delta)^{\frac{\alpha}{2}}u(x):=\mathcal{F}^{-1}\left(1+\lvert \xi \rvert^{\alpha}\hat{u}(\xi)\right)(x),$$
where $\hat{\xi}=\mathcal{F}(u(x))$ is the Fourier transformation of function $u(x)$. Also, we will use the following formula as the one given in \cite{caffarelli2007extension}:
$$-(-\Delta)^{\frac{\alpha}{2}}u(x)=C_{\alpha,n}P.V.\int_{\mathbb{R}^{n}}\frac{u(y)-u(x)}{\lvert y-x\rvert^{n+\alpha}}dy,$$
where $C_{\alpha,n}=\frac{2^{\alpha-1}\alpha\Gamma((n+\alpha)/2)}{\pi^{2/d}\Gamma(1-\alpha/2)}$ is normalization constant P.V. denotes the Cauchy principle value.
The unknown function $v(x,t)$ is given by the fundamental solution
\begin{equation}
    v(x,t)=\begin{cases}
    \frac{1}{n(n-2)\omega(n)}\int_{\mathbb{R}^{n}}\frac{1}{\lvert x-y\rvert^{n-2}}u(y)dy, &{ \rm if} \ n\geq 3,\\
    -\frac{1}{2\pi}\int_{\mathbb{R}^{n}}\ln\lvert x-y \rvert u(y,t)dx, & {\rm if} \ n=2,
    \end{cases}
\end{equation}
where $\omega(n)=\frac{\pi^{\frac{n}{2}}}{\Gamma(n/2+1)}$ is the wolume of the $n$-dimensional unit ball. 

Before we go to further results, we recast $v$ in (\ref{1.1}) as $v=\Psi*u$ where $\Psi(x)$ is the Newton potential, and it can be represented as
\begin{equation}
    \Psi(x)=\begin{cases}
    \frac{C_{n}}{\lvert x\rvert^{n-2}}, &{ \rm if} \ n\geq 3,\\
    -\frac{1}{2\pi}\ln\lvert x\rvert , & {\rm if} \ n=2.
    \end{cases}
\end{equation}
Thus we have the attractive force
\begin{equation}
    F(x)=\nabla \Psi(x)=-\frac{C_{*}x}{\lvert x\rvert^{n}},\ \ \forall x\in \mathbb{R}^{n}\ {0},
\end{equation}
where $C_{*}=\frac{\Gamma(n/2)}{2\pi^{2/n}}$. Moreover $\nabla v=F*u$.

Fractional derivatives are differential operators that extend the definition of integer derivatives to fractional derivatives, which are used to describe nonlocal effects in time and space \cite{gorenflo2015time,kilbas2006theory}. In recent years, there has been a lot of research on fractional calculus in many application fields, such as physics, biomedical engineering, chemistry, and control theory \cite{klafter2005anomalous,li2017fractional,lorenzi1988inverse,delbosco1996existence}.The most important advantage of introducing fractional differential equations in these applications is their non-locality, which means that the next state of a system depends not only on its current state, but also on all its historical states. Since the behavior of most biological systems has memory and after-effects. Therefore, it is more advantageous to model these systems by fractional differential equations. The fractional derivative in time pass is used to simulate memory effects that can be seen everywhere \cite{caputo1967linear,li2018generalized}. When the diffusion velocity of the particle plume is inconsistent with the diffusion velocity of the Brownian motion model, both the time fractional derivative and the spatial fractional derivative can be used for anomalous diffusion or dispersion. The spatial fractional derivative term in the diffusion equation can describe the macroscopic phenomenon of transport and usually leads to the occurrence of superdiffusion.

In the last few decades, a lot of work has been done to study fractional solutions to partial differential equations\cite{eidelman2004cauchy,kemppainen2017representation,vazquez2014recent}. For example, Shu et al.\cite{shu2011existence} studied the existence of mild solutions for a class of impulse fractional semi-linear partial differential equations.zacher \cite{zacher2013giorgi} studied the regularity of weak solutions of linear diffusion equations in Riemann-Liouville time fractional derivatives in bounded regions in $\mathbb{R}^{n}$ with arbitrary time-order $\alpha$-divergence forms. Allen, Caffarelli and Vasseur \cite{allen2017porous} discussed parabolic problems with Caputo-type fractional time derivatives. At present, the Keller-Segel equation with time fractional derivative and fractional Laplace term also has many conclusions. Huang and Liu \cite{huang2016well} used the basic energy method to obtain a priori estimate, then constructs the regularization equation, and obtains the existence of the regularized smooth solution, and finally uses the Aubin-Lions-Dubinskii lemma to prove that the Keller-Segel equation containing the nonlocal diffusion term $-\nu(-\Delta)^{\frac{\alpha}{2}}u$ is weak The existence of the solution and the uniqueness and stability of the weak solution under the Wasserstein metric are also proved by linking the KS equation to a self-consistent stochastic process driven by a rotationally invariant $\alpha$-stable Levy process $L_{\alpha}(t)$. \cite{biler2009two} constructed the global in time solutions of fractional diffusion $1<\alpha \leq 2$, and analyzed that all solutions are globally bounded in time in $n=1$ dimension space. Zhou et al. \cite{zhou2018weakness} studied the time fractional Keller-Segel model, and obtained the existence of weak solutions using the Faedo-Galerkin approximation method combined with the compactness criterion. \cite{biler2010blowup} studied the application of linear analytic semigroup theory in the process of quasi-linear evolution, and discussed the initial value $u_0$ in the critical Besov space $\dot{B}_{2,r}^{1-\alpha}(\mathbb{R}^{2})$ for $r\in [1,\infty)$ and $1<\alpha<2$. Using the fixed point theorem to study the existence of solutions to time-fractional KS models with Caputo-Fabrizlo fractional derivatives can be found in \cite{atangana2015analysis}. \cite{kumar2017new} considered the dual time fractional KS equation of $-1<\beta\leq n$, and used the modified homotopy analysis transformation method to solve the one-dimensional time fractional Keller-Segel model. In two-dimensional space, Li and Liu \cite{li2018some} obtained a compactness criterion similar to the Aubin-Lions lemma, which can be used to prove the existence of weak solutions to time fractional partial differential equations, and also proved the existence of weak solutions to time fractional Keller-Segel equations on R as The specific application of the compactness criterion. Li et al. \cite{li2018cauchy} discussed the Cauchy problem of the Keller-Segel type diffusion equation extended in space-time fractional order, obtained the local and global existence and uniqueness of the mild solution by using the $L^{r}-L^{q}$ estimation of the fundamental solution and the fixed point theorem, and also studied the mass conservation, non-negative Sexuality and demolition issues.

Model (\ref{1.1}) can be seen as a generalization of the classical parabolic-elliptic Keller-Segel model.
\begin{equation}\label{1.2}
\left\{
\begin{aligned}
    & u_{t} = \Delta u - \nabla \cdot (u\nabla v),  & x\in \mathbb{R}^{n}, t>0,\\
    & 0=\Delta v+u, & x\in \mathbb{R}^{n}, t>0,\\
    & u(x,0)  = u_0,  & x\in \mathbb{R}^{n},\\
\end{aligned}
\right.
\end{equation}
Mathematical models related to the theory of biochemotaxis can be traced back to \cite{patlak1953random} in 1950. The original parabolic-parabolic Keller-Segel model was introduced by Keller and Segel \cite{keller1970initiation} in 1970 to describe the chemotactic migration process of cell slime molds, while the parabolic-ellipse model was proposed by J$\ddot{a}$gger and Luckhaus \cite{jager1992explosions} in 1992 due to the very large diffusion coefficient of chemical substances. The second equation of (\ref{1.2}) simulates the fact that the cell slime mold produces chemicals, and the term 
$- \nabla \cdot (u\nabla v)$ comes from the chemotactic migration process, that is, from the bacteria moving with the chemical gradient. Therefore, $- \nabla \cdot (u\nabla v)$ is the aggregate term, which cannot be bounded in all cases, and in some cases the buckling state causes the solution to blow up in finite time \cite{winkler2013finite}.

J$\ddot{a}$gger and Luckhaus \cite{jager1992explosions} first construct the burst solution of (\ref{1.2}) and propose a conjecture that the size of u will affect whether the solution bursts or exists globally in finite time. \cite{nagai2001blowup} verified this conjecture. In the two-dimensional space, if $\int_{B_{L}}u_{0}dx<8\pi$, the global existence of the solution of (\ref{1.2}) in time is obtained. If $u_0$ satisfies $\int_{\Omega}u_0dx>8\pi$ and $\int_{\Omega}u_0\lvert x-q \rvert^{2}dx$ is small for $q\in\Omega$, the solution of (\ref{1.2}) will blow-up in finite time. It is well known that the properties of a model solution depend on the spatial dimension of the solution, see \cite{childress1984chemotactic}. The work in \cite{horstmann20031970,horstmann20041970} gives the detailed content, that is, $n=1$, the smooth solution exists globally; $n=2$, the existence of the solution depends on the size of the initial mass $\Vert u_{0}\Vert_{L^{1}(\mathbb{R}^{n})}$, and the solution exists globally for small initial mass, and the large initial mass exists. The initial mass solution will explode; $n=3$; the existence or non-existence of the solution depends on the size of $\Vert u_{0}\Vert_{L^{\frac{n}{2}}(\mathbb{R}^{n})}$. On the other hand, the study of the overall existence of the solution of the Keller-Segel model and the optimal initial critical mass for blasting is also an important aspect, which is also related to the spatial dimension. For example, \cite{osaki2001finite} proved the global existence of weak solutions in one-dimensional space, and in two-dimensional space, there is a threshold that determines whether the solution exists in totality or blows up in time. $n=3$, no such threshold exists for  work in \cite{winkler2010aggregation}. This problem has been extensively studied in two-dimensional space, and it has been proved that $8\pi$ is the optimal initial critical mass for the Keller-Segel model in two-dimensional space, see \cite{blanchet2008infinite,blanchet2006two}. It is given in \cite{dolbeault2004optimal}: if the initial mass $M_0<8\pi$, the weak solution exists as a whole; if $M_0>8\pi$, the solution blows up in a finite time. But there is less discussion about the optimal initial criticality for high-dimensional cases. Regarding other properties of the solution of the model (\ref{1.2}),  \cite{liu2016refined,liu2016note} proved the uniqueness and $L^{\infty}$-bound of the weak solution .

Our model (\ref{1.1}) replaces the time derivative in the classical model (\ref{1.2}) with the Caputo fractional derivative, replaces the Laplacian with the fractional Laplacian, takes into account memory effects and non-locality, and describes chemotaxis of biological systems for paradoxical diffusion and memory effects. In addition, to the first equation of (\ref{1.2}), a characterization-logistic source term, which is very important for the species to live or die, is added. Inspired by the above literature, we study the existence, $L^{\infty}$-bound and uniqueness of weak solutions of model (\ref{1.1}), and we also establish the blow-up criterion for weak solutions. The existence of weak solutions is mainly discussed in terms of the damping coefficient $b$. We address the global existence of weak solutions of fractional Keller-Segel system (\ref{1.1}) in a standard approach. We start with the construction of some a prior estimates of the solutions to the mollified Equations \ref{2.13}.

Suppose function $h(x,t)$ is smooth and all its derivatives are bounded. To study the existence of weak solutions, we first of all investigate the following fractional advection diffusion equation:
\begin{equation}\label{2.13}
\left\{
\begin{aligned}
    &  \partial_{t}^{\beta} u+(-\Delta)^{-\frac{\alpha}{2}}u + \nabla \cdot (uh(x,t))-F(x,t)= 0, \\
    & u(x,0)  = u_0(x).
\end{aligned}
\right.
\end{equation}
From \cite{taylorremarks}, taking Laplace transform on both of the first equation of (\ref{2.13}), one has the following analogy of Duhamel's principle
\begin{equation}\label{2.14}
    u(x,t)=E_{\beta}(-t^{\beta}A)u_{0}+\beta\int_{0}^{t}s^{\beta-1}E_{\beta}^{'}(-t^{\beta}A)(-\nabla \cdot (ua)|_{t-s}+F|_{t-s})ds.
\end{equation}
The $L^{\infty}$-bound of the weak solution is mainly obtained by establishing the fractional differential inequality and then using the Moser iteration method. Using hypercompression estimation, the conclusion that weak solutions exist must be unique.

We note that, to the best of our knowledge, the properties of the existence $L^{\infty}$-bound and uniqueness of weak solutions to spatiotemporal fractional Keller-Segel type diffusion equations with logistic source terms have not been discussed so far. Therefore, in this paper, we are interested in making strict requirements on the global existence of weak solutions and their $L^{\infty}$-norm uniform boundedness and uniqueness. The structure of this paper is as follows: Section 2 presents some basic definitions and properties, which will be used in the subsequent proof process, and Section 3 is devoted to discussing the global existence of weak solutions to model (\ref{1.1}). At the same time, we also prove the  $L^{\infty}$-bound of weak solutions. In Section 4, we use the hyper-contractivity of solutions to prove the uniqueness. In Section 5, we also establish the blow-up criterion for weak solutions. The main results of this paper are as follows:

\begin{thm}\label{1.1.1}
Assume that the initial data $u_{0} \in L^{1}(\mathbb{R}^{n}) \cap L^{\infty}(\mathbb{R}^{n})$, and the parameters are satisfied respectively
    \begin{enumerate}[(i)]
        \item $a \geq 0, b \geq 1$
        \item $a \geq 0, 1-\frac{\alpha}{n} < b < 1$
    \end{enumerate}
     Then there exists a non-negative global weak solution $(u,v)$ of (\ref{1.1}). and there is a constant $C$ such that
    \begin{equation}
         \Vert u \Vert_{L^{\infty}(0,T;L^{\infty}(\mathbb{R}^{n}))} \leq C.
    \end{equation} 
\end{thm}

\begin{thm}\label{1.1.2}
 Assume that $u_0\in L^{1}(\mathbb{R}^{n}) \cap L^{\frac{n}{\alpha}}(\mathbb{R}^{n})$ and $\Vert u_0 \Vert_{L^{\frac{n}{\alpha}}}\leq \frac{{C_{*}}}{2}$, where $C_{*}=\frac{2 S_{\alpha,n}^{-2}\left(\frac{n}{\alpha}-1\right)}{\frac{n}{\alpha}\left(\frac{n}{\alpha}-1-\frac{n}{\alpha}b\right)}$, the parameter $b$ satisfies $0< b \leq 1-\frac{\alpha}{n}$. For any fixed $T>0$, let $a$ satisfy $0 \leq a < \left(\frac{(\frac{C_{*}}{2})^{\frac{n}{\alpha}}\Vert u_0\Vert_{L^{2}(\mathbb{R}^{n})}^{\frac{n}{\alpha}}}{C_{0}T}\right)^{\tau_5}$, where $\tau_5=\frac{\alpha^2}{n^2-\alpha(n-\alpha)(1+\sigma)}$. Then the system (\ref{1.1}) has a global weak solution.
\end{thm}

\begin{thm}\label{1.1.3}
Suppose the initial data $u_0 \in L^{1}(\mathbb{R}^{n}) \cap L^{\infty}(\mathbb{R}^{n})$, $u$ is a weak solution of the system (\ref{1.1}) under the condition of Theorem \ref{1.1.1}, then the weak solution $u$ is unique in its corresponding existence space.  
\end{thm}

\section{Notations and preliminaries}
The purpose of this section is to introduce some concepts and list some of the results used in the paper.

$H^{1,1}(0,T)=\{\rho(t)\in L^{1}(0,T):\dot{\rho}(t)\in L^{1}(0,T)\}$, where $\dot{\rho}(t)=\frac{\partial}{\partial t}\rho(t)$. For $s\in \mathbb{R}$, $1<p<\infty$, non-integer power Sobolev space $W^{s,p}(\mathbb{R}^{n})$ is defined as $W^{s,p}(\mathbb{R}^{n})=\left\{u\in L^{p}(\mathbb{R}^{n}): \mathcal{F}^{-1}\left((1+\lvert \xi \rvert^{2})^{\frac{s}{2}}\hat{u}(\xi)\right)\in L^{p}(\mathbb{R}^{n})\right\}$ with the norm
$$\Vert f \Vert_{W^{s,p}(\mathbb{R}^{n})}=\left\Vert \mathcal{F}^{-1}\left((1+\lvert \xi \rvert^{2})^{\frac{s}{2}}\hat{f}(\xi)\right) \right\Vert_{L^{p}(\mathbb{R}^{n})}.$$ 
The space $C([0,T];H^{\gamma}(\mathbb{R}^{n}))$ comprises all continuous functions $u: [0,T]\to H^{\gamma}(\mathbb{R}^{n})$ with
$$\Vert u \Vert_{C([0,T];H^{\gamma}(\mathbb{R}^{n}))}:=\underset{0\leq t \leq T}{\sup}\Vert u \Vert_{H^{\gamma}(\mathbb{R}^{n})}<\infty.$$
In addition, denote $A:=(-\Delta)^{\frac{\alpha}{2}}$ to be the fractional Laplacian. Throughout this work, the constant $C$ denotes generic positive constant in the later proofs which may vary from line to line.

Let us first recall the following definition of limit.
\begin{definition}\cite{li2018some}\label{2.1.1}
    Let $B$ be a space. For a function $u\in L_{loc}^{1}(0,T;B)$, if there exists $u_0\in B$ such that
    $$\underset{t\to 0^{+}}{\lim}\frac{1}{t}\int_{0}^{t}\Vert u(s)-u_0\Vert_{B}ds=0,$$
we call $u_0$ the right limit of $u$ at $t=0$, denoted by $u(0+)=u_0$.

we define $u(T-)$ to be the constant $u_{T}\in B$ such that
$$\underset{t\to T^{-}}{\lim}\frac{1}{T-t}\int_{t}^{T}\Vert u(s)-u_T\Vert_{B}ds=0,$$
For $\beta >-1$, as dicussed in \cite{li2018some}, we define $\{g_{\beta}\}$ as the convolution kernels
\begin{equation}
    g_{\beta}(t):=\begin{cases}
    \frac{\theta(t)}{\Gamma(\beta)}t^{\beta-1}, &\beta>0;\\
    \delta(t), &\beta=0;\\
    \frac{1}{\Gamma(1+\beta)}D(\theta(t)t^{\beta}), &\beta\in(-1,0),
    \end{cases}
\end{equation}
where $\theta(t)$ is the standard Heaviside step function and $D$ represents the distributional derivative. $g_{\beta}$ can also be defined for $\beta\leq -1$ and consequently
$$g_{\beta_{1}}*g_{\beta_{2}}=g_{\beta_1+\beta_2}, \ \ \forall \beta_{1},\beta_{2}\in \mathbb{R}.$$
Correspondingly, the time-reflected group:
$$\tilde{\Psi}:=\{\tilde{g}_{\alpha}:\tilde{g}_{\alpha}=g_{\alpha}(-t), \alpha\in \mathbb{R}\}.$$
Clearly, ${\rm supp\ }\tilde{g}\subset (-\infty,0]$ and for $\gamma\in(0,1)$ the following equality is true:
$$\tilde{g}_{\gamma}(t)=-\frac{1}{\Gamma(1-\gamma)}D(\theta(-t)(-t)^{\gamma})=-D\tilde{g}_{1-\gamma}(t).$$
\end{definition}

\begin{definition}\cite{li2018some}
    Let $0<\beta<1$. Consider $u\in L_{loc}^{1}(0,T;\mathbb{R})$ such that $u$ has a right limit $u(0+)$ at $t=0$ in the sense of Definition \ref{2.1.1}. The $\beta$th order Caputo derivative of $u$, a distribution in $\mathscr{D}'(\mathbb{R})$ with support in $[0,T)$, is defined by
    $$\partial_{c}^{\beta}u:=J_{-\beta}u-u_0g_{1-\beta}=g_{-\beta}*\left(\theta(t)(u-u_0)\right),$$
    where $J_{-\beta}$ denotes the fractional intergral operator
    $$J_{\beta}u(t)=\frac{1}{\Gamma(-\beta)}\int_{0}^{t}(t-s)^{\alpha-1}u(s)ds.$$
    
    Similarly, the $\beta$th order right Caputo derivative of $u$ is a distribution in $\mathscr{D}'(\mathbb{R})$  with support in $(-\infty,T]$, given by
    $$\tilde{\partial}_{T}^{\beta}u:=\tilde{g}_{-\beta}*\left(\theta(T-t)(u(t)-u(T-))\right).$$
\end{definition}

\begin{definition}\cite{li2018some}
     Let $0<\beta<1$. Consider $u\in L_{loc}^{1}(0,T;\mathbb{R})$ such that $u$ has a right limit $u(0+)$ at $t=0$ in the sense of Definition \ref{2.1.1}. The $\beta$th order right Caputo derivative of $u$, a distribution in $\mathscr{D}'(\mathbb{R})$ with support in $[0,T)$, given by
     \begin{equation}
         \tilde{\partial}_{T}^{\beta}u:=\tilde{g}_{-\beta}*\left(\theta(T-t)(u(t)-u(T-))\right).
     \end{equation}
\end{definition}

\begin{definition}\cite{li2018some}
    Let $B$ be a Banach space and $u\in L_{loc}^{1}(0,T;B)$. Let $u_0\in B$. We define the weak Caputo derivative of $u$ associated with initial data $u_0$ to be $\partial_{t}^{\beta}u\in  \mathscr{D}'.$ such that for any test function $\phi\in C_{0}^{\infty}((-\infty,T);\mathbb{R})$
    $$\langle \partial_{t}^{\beta}u, \phi\rangle:=\int_{0}^{T}(u-u_0)\theta(t)\left(\tilde{\partial}_{T}^{\beta}\phi\right)dt=\int_{0}^{T}(u-u_0)\tilde{\partial}_{T}^{\beta}\phi dt,$$
    where $\mathscr{D}'(\mathbb{R})=\left\{\nu|\nu:C_{c}^{\infty}((-\infty,T);\mathbb{R})\to B \ is\ a \ bound\ linear\  operator\right\}.$ We call the weak Caputo derivative $\partial_{t}^{\beta}u$ associated with initial value $U_0$ the Caputo derivative of $u$ if $u(0+)=u_0$ in the sense of Definition \ref{2.1.1} under the norm of the underlyiing Banach space $B$.
\end{definition}
\begin{definition}\cite{huang2016well}(Weak Solution)\label{2.1}
Let $0\leq u_0(x) \in L^{1}(\mathbb{R}^{n})\cap L^{q}(\mathbb{R}^{n})$ be the initial data , where $q>1$ and $T \in (0,\infty)$. We say $u(x.t)$ is a weak solution to (\ref{1.1}) with initial data $u_{0}(x)$ if it satisfies
\begin{enumerate}[(1)]
\item for any $T>0$, $q_1>1$
\begin{equation}\label{2.2}
    u(x,t)\in L^{\infty}(0,T;L^{1}(\mathbb{R}^{n})\cap L^{q}(\mathbb{R}^{n})),\ \ u\in L^{2}(0,T:H^{\frac{\alpha}{2}}(\mathbb{R}^{n})),
\end{equation}
\begin{equation}\label{2.3}
    \partial_{t}^{\beta}u\in L^{q_1}(0,T;W^{-\alpha,p_{1}}(\mathbb{R}^{n}))\ \ \ {\rm for}\ {\rm some}\  p_{1}\geq 1.
\end{equation}
\item For any $\phi(x) \in C_{c}^{\infty}(\mathbb{R}^{n})$ and $0< t < T$, it holds
\begin{align}\label{2.4}
        &\int_{0}^{T}\int_{\mathbb{R}^{n}}(u-u_0)\phi(x,t)\tilde{\partial}_{T}^{\beta}\phi(x,t)dxdt\notag\\
        ={}& - \int_{0}^{T}\int_{\mathbb{R}^{n}}\left[u(x,t)D^{\alpha}\phi(x)\right]dxdt+\int_{0}^{T}\int_{\mathbb{R}^{n}}u(x) \nabla v(x)\cdot \nabla \phi(x)dxdt\notag\\
        +{} & a\int_{0}^{T}\int_{\mathbb{R}^{n}}u(x)\phi(x)dxdt-b\int_{0}^{T}\int_{\mathbb{R}^{n}}u^{2}(x)\phi(x)dxdt.
\end{align}
\item $v(x,t)$ is the chemicial substance concentration associated with $u$ and given by
\begin{equation}
    \nabla v=\int_{\mathbb{R}^{n}}F(x-y)u(y,t)dy.
\end{equation}
\end{enumerate}
\end{definition}

\begin{definition}\cite{li2018cauchy}\label{2.3.1}
  Let $X$ be a Banach space over space and time. We call $u\in X$ is a mild solution to 
(\ref{2.13}) if $u$ satisfies the integral equation (\ref{2.14}) in $X$.  
\end{definition}

\begin{definition}\cite{huang2016well}(Weak Solution)\label{2.4.1}
Let $0\leq u_0(x) \in L^{1}(\mathbb{R}^{n})\cap L^{\frac{n}{\alpha}}(\mathbb{R}^{n})$ be the initial data and $T \in (0,\infty)$. We say $u(x.t)$ is a weak solution to (\ref{1.1}) with initial data $u_{0}(x)$ if it satisfies
\begin{enumerate}[(1)]
\item for any $T>0$, $q_1>1$
\begin{equation}\label{2.4.2}
    u(x,t)\in L^{\infty}(0,T;L^{1}(\mathbb{R}^{n})\cap L^{\frac{n}{\alpha}}(\mathbb{R}^{n})),\ \ u\in L^{2}(0,T:H^{\frac{\alpha}{2}}(\mathbb{R}^{n})),
\end{equation}
\begin{equation}\label{2.4.3}
    \partial_{t}^{\beta}u\in L^{q_1}(0,T;W^{-\alpha,r_2}(U))\ \ \ r_2=\min\left\{2,\frac{n(n+\alpha)}{n(\alpha+1)+\alpha}\right\}.
\end{equation}
\item For any $\phi(x) \in C_{c}^{\infty}(\mathbb{R}^{n})$ and $0< t < T$, (\ref{2.4}) holds.

\end{enumerate}
\end{definition}

\begin{definition}\cite{zhou2018weakness}\label{2.5}
    Assume that $X$ is a Banach space and let $u: [0,T] \to X$. The Riemann–Lioville fractional derivatives of $u$ for order $\beta \in \mathbb{C} (Re(\beta)>0)$ are defined by
    \begin{gather*}
        _{0}\textrm{D}_{t}^{\beta}u(t)=\frac{1}{\Gamma(1-\beta)}\frac{d}{dt}\int_{0}^{t}(t-s)^{-\beta}u(s)ds\\
         _{0}\textrm{D}_{T}^{\beta}u(t)=\frac{-1}{\Gamma(1-\beta)}\frac{d}{dt}\int_{t}^{T}(t-s)^{-\beta}u(s)ds,\ t>0.
    \end{gather*}
    where $\Gamma(1-\beta)$ is the Gamma function. The above integrals are called the left-sided and the right-sided the Riemann-Liouville fractional derivatives.
\end{definition}

\begin{definition}\cite{zhou2018weakness}\label{2.6}
    (i) Assume that $X$ is a Banach space and let $u: [0,T] \to X$. The Caputo fractional derivative operators of $u$ for order $\beta \in \mathbb{C} (Re(\beta)>0)$ are defined by
     \begin{gather*}
        _{0}^{C}\textrm{D}_{t}^{\beta}u(t)=\frac{1}{\Gamma(1-\beta)}\int_{0}^{t}(t-s)^{-\beta}\frac{d}{dt}u(s)ds\\
         _{0}^{C}\textrm{D}_{T}^{\beta}u(t)=\frac{-1}{\Gamma(1-\beta)}\int_{t}^{T}(t-s)^{-\beta}\frac{d}{dt}u(s)ds,
    \end{gather*}
    where $\Gamma(1-\beta)$ is the Gamma function.  The above integrals are called the left-sided and the right-sided the Caputo fractional derivatives.\\
    (ii) For $u:[0,\infty) \times \mathbb{R}^{n} \to \mathbb{R}$, the left Caputo fractional derivative with respect to time $t$ of $v$ is defined by
    $$\partial_{t}^{\beta}u = \frac{1}{\Gamma(1-\beta)}\int_{0}^{t}(t-s)^{-\beta}\frac{\partial}{\partial s}u(s,x)ds,\ \ t>0.$$
\end{definition}

\begin{definition} \cite{YU2018306}\label{2.7}
    The one- and two-parameter Mittag-Leffler functions can be  defined by 
    \begin{equation}
    \left\{
        \begin{array}{lr}
         E_{\beta}(z) = \sum_{k=0}^{\infty}\frac{z^{k}}{\Gamma(\beta k+1)}, & \beta>0.\\
         E_{\beta,\gamma}(z) = \sum_{k=0}^{\infty}\frac{z^{k}}{\Gamma(\beta k+\gamma)}, & \beta,\gamma >0.
         \end{array}
    \right.
    \end{equation}
    we will dennote $E_{\beta}(z)\overset{\Delta}{=}E_{\beta,1}(z)$.
\end{definition}

\begin{lemma}\cite{li2018some}\label{2.8.1}
Let $0<\beta<1$. If the mapping $u:[0,T)\to B$ satisfies $u\in C^{1}(0,T;B)\cap C(0,T;B)$ and $u \mapsto E(u)\in \mathbb{R}$ is a $C^{1}$ convex functional on $B$, then
$$\partial_{t}^{\beta}u(t)=\frac{1}{\Gamma(1-\beta)}\left(\frac{u(t)-u(0)}{t^{\beta}}+\beta\int_{0}^{t}\frac{u(t)-u(s)}{(t-s)^{\beta+1}}ds\right),$$
and 
$$\partial_{t}^{\beta}E(u(t))\leq \langle D_{u}E(u(t)), \partial_{t}^{\beta}u\rangle,$$
where $D_{u}E(\cdot): B\to B'$ is Fr\'{e}chet and $\langle \cdot, \cdot\rangle$ is understood as the dual pairing between $B'$ and $B$.
\end{lemma}

\begin{lemma}\cite{li2018cauchy}\label{2.8.2}
    Let $\beta\in(0,1)$.\\
(i) Let $u \in L_{loc}^{1}(0,T;\mathbb{R})$. Assume that weak Caputo derivative for an assigned initial value $u_{0}\in \mathbb{R}$ is $\partial_{t}^{\beta}u$. As linear functions on $C_{c}^{\infty}(-\infty,T;\mathbb{R})$, we have
\begin{equation}
    (u-u_0)\theta(t)=g_{\beta}*(\partial_{t}^{\beta}u).
\end{equation}
(ii) If $f(t):=(\partial_{t}^{\beta}u)\in L_{loc}^{1}(0,T;B)$, then 
\begin{equation}
    u(t)=u_{0}+\frac{1}{\Gamma(\gamma)}\int_{0}^{t}(t-s)^{\gamma-1}f(s)ds, \ a.e.\ on (0,T),
\end{equation}
where the integral is understood as the Lebesgue integral.
\end{lemma}

\begin{lemma}\cite{li2018some}\label{2.8.3}
Let $0<\beta<1$, $T>0$. Assume that $z(t)\in C([0,T];\mathbb{R})$. Suppose $f(t,x)$ is a continuous function, locally Lipschitz in $x$, such that $\forall t\geq 0$, $x\leq y$ implies $f(t,x)\leq f(t,y)$. If $f(t,z)-(\partial_{t}^{\beta}v)$ is a nonnegative distribution, then $v\leq u$ for $t\in[0,\min(T,T_{b}))$, where $u$ is the solution to the ODE
\begin{equation}
    \partial_{t}^{\beta}u=f(t,z),\ u(0)=z(0)
\end{equation}
and $T_{b}$ is the largest existence time for $u$.
\end{lemma}

\begin{lemma}\cite{bonforte2014quantitative}(Sobolev inequality)\label{2.8}
Let $0<\frac{\alpha}{2}\leq 1$ and $\alpha <n$. Then
\begin{equation}
    \Vert f \Vert_{\frac{2n}{n-\alpha}}\leq S_{\alpha,n}\Vert D^{\frac{\alpha}{2}}f \Vert_{2}
\end{equation}
where the best constant is given by
$$S_{\alpha,n}^{2}:=2^{-\alpha}\pi^{\frac{\alpha}{2}}\frac{\Gamma(\frac{n\alpha}{2})}{\Gamma(\frac{n+\alpha}{2})}\left[\frac{\Gamma(n)}{\Gamma(\frac{n}{2})}\right]^{\frac{\alpha}{n}}=\frac{\Gamma(\frac{n-\alpha}{2})}{\Gamma(\frac{n+\alpha}{2})}\lvert \mathbb{S}^{n-1}\rvert^{-\frac{\alpha}{n}}.$$
\end{lemma}

\begin{lemma}\cite{bonforte2014quantitative}(Stroock-Varopoulos’ inequality)\label{2.9}
Let $0<\frac{\alpha}{2}<1$, $p>1$, then
\begin{equation}\label{2.9.1}
    -\int_{\mathbb{R}^{n}}\lvert f\rvert^{p-2}fD^{\alpha}fdx\leq -\frac{4(p-1)}{p^{2}}\Vert D^{\frac{\alpha}{2}}f^{\frac{p}{2}}\Vert_{2}^{2}
\end{equation}
for all $f\in L^{p}(\mathbb{R}^{n})$ such that $D^{\alpha} f\in L^{p}(\mathbb{R}^{n})$.
\end{lemma}

\begin{lemma}\cite{alsaedi2015maximum}\label{2.10}
For any function $\tau(t)$ absolutely continuous on $[0,T]$, one has the inequality
    $$\tau^{p-1}(t)\partial_{t}^{\beta}\tau(t)\geq \frac{1}{p}\partial_{t}^{\beta}\tau^{p}(t).$$ 
\end{lemma}

\begin{lemma} \cite{alikhanov2010priori}\label{2.11}
     Let a non-negative absolutely continuous functon $z(t)$ satisfy the inequality
    \begin{equation}
        \partial_{t}^{\beta}z(t) \leq c_1z(t)+c_2(t),\ \ \ \ 0 < \beta \leq 1
    \end{equation}
    for almost all $t$ in $[0,T]$, where $c_1 >0$ and $c_2(t)$ is an integrable non-negative function on $[0,T]$. Then
    \begin{equation}
        z(t) \leq z(0)E_{\beta}(c_1t^{\beta})+\Gamma(\beta)E_{\beta.\beta}(c_1t^{\beta})D_{t}^{-\beta}c_2(t)
    \end{equation}
    where $E_{\beta}(y)$ and $E_{\beta,\gamma}(y)$ are the Mittag-Leffler functions.
\end{lemma}

\begin{lemma} \cite{zhou2018weakness}\label{2.12}
    For $u(t) \geq 0$,
    $$_{0}^{c}\textrm{D}_{t}^{\beta}u(t)+c_1u(t)\leq c_2(t)$$
    for almost all $t\in[0,T]$, where $c_1 >0$, and the function $c_2(t)$ is non-negative and integrable for $t\in[0,T]$. Then
    $$u(t)\leq u(0)+\frac{1}{\Gamma(\beta)}\int_{0}^{t}(t-s)^{\beta-1}c_2(s)ds.$$
\end{lemma}

\section{The existence of weak solutions and the \texorpdfstring {$L^{\infty}$}{lg}-bound}
This section mainly proves the existence and $L^{\infty}$-bound of weak solutions of the model (\ref{1.1}). We divide the discussion into three cases: (i) $a\geq 0, b\geq 1$, the initial value is arbitrary; (ii) $a\geq 0$, $1-\frac{\alpha}{n} < b <1$ initial value is arbitrary; (iii) $0<b\leq 1-\frac{\alpha}{n}$, small initial value, satisfying $\Vert u_{0}\Vert_{L^{\frac{n}{\alpha}}} < C_{*}$ , and the birth rate is small.

We intend to discuss the problem that if the fact that the initial value $U_0$ of (\ref{2.13}) is non-negative can imply that $u(x,t)$ remains non-negative for every $0<t<T$. 

Given a function $u$, we define $u^{-}=-\min(u,0)\geq 0$, $u^{+}=\max(u,0)$ so that
$$u=u^{+}-u^{-}.$$

\begin{lemma}\cite{li2018cauchy}
Let $A=(-\Delta)^{\frac{\alpha}{2}}$. Suppose $u\in L^{2}(\mathbb{R}^{n})$, then for any $t>0$,
\begin{equation}
    \langle e^{-tA}u, u^{+}\rangle \leq \Vert u^{+} \Vert_{L^{2}}^{2}=\int_{\mathbb{R}^{+}}uu^{+}dx.
\end{equation}
\end{lemma}

\begin{lemma}\cite{li2018cauchy}
Let $A=(-\Delta)^{\frac{\alpha}{2}}$. If $u\in H^{\delta, 2}$ for $\delta\in [0,1]$ and $Au\in H^{-\delta,2}(\mathbb{R}^{n})$, then we have 
\begin{equation}
    \langle Au, u^{+}\rangle \geq 0,\ \ \langle Au, u^{-}\rangle \leq 0.
\end{equation}
\end{lemma}

\begin{lemma}\cite{li2018some}\label{3.1.1.1.1}
Let $T>0$, $0<\beta<1$ and $1<p<\infty$. Let $M,B$ and $Y$ be Banach spaces. $M\hookrightarrow B$ compactly and $B\hookrightarrow Y$ continuously. Suppose that $M\subset L_{loc}^{1}((0,T);M)$ satisfies the following conditions,\\
(i) here exists $C_1>0$ such that for any $u\in W$,
$$\underset{t\in(0,T)}{\sup}J_{\beta}(\Vert u \Vert_{M}^{p})=\underset{t\in(0,T)}{\sup}\int_{0}^{t}(t-s)^{\beta-1}\Vert u \Vert_{M}^{p}(s)ds\leq C_1;$$
(ii) there exists $r\in \left(\frac{p}{1+p\beta},\infty\right)\cap [1,\infty)$ and $C_2>0$ such that for any $u\in W$, there is an assignment of initial value to make the weak Caputo derivative satisfy
$$\Vert \partial_{t}^{\beta}u\Vert_{L^{r((0,T);Y)}}\leq C_2.$$
Then $W$ is relatively compact in $L^{p}((0,T);B)$.
\end{lemma}

\begin{lemma}\label{3.1.1.1.1.1}
Suppose $n\geq 2$, $0<\beta<1$, $1<\alpha < 2$. Let $p\in (\frac{\alpha}{n},\infty)\cap[\frac{2n}{n+1},\frac{n}{2})$, for $u_{0}\in L^{p}(\mathbb{R}^{n})\cap H^{\gamma}(\mathbb{R}^{n})$. Let $(u,v)$ be a mild solution of (\ref{2.13}) in $C([0,T];L^{p}(\mathbb{R}^{n}))$ with initial data $u_0$ in sense of Definition \ref{2.3.1}. If we also have $u_0\geq 0$ for all $t$ in the interval of existence, we have 
\begin{equation}
    u(x,t)\geq 0.
\end{equation}
\end{lemma}
\begin{remark}
We prove the non-negativity of the mild solution of the nonlinear fractional convection-diffusion equation, which plays an important role in the proof of the latter lemma.
\end{remark}
\begin{pro}
We introduce a mollifier $J_{\varepsilon}(x)=\frac{1}{\varepsilon^{n}}J(\frac{x}{\varepsilon})$ and consider $-\Delta v_{\varepsilon}=u_{\varepsilon}*J_{\varepsilon}$. We adopt the method in \cite{li2018cauchy}, 
$$T_{b}=\sup\{T>0: (\ref{1.1})\ {\rm has}\ {\rm a}\ {\rm unique}\ {\rm mild}\ {\rm solution}\ {\rm in}\  C[0,T;L^{p}(\mathbb{R}^{n})]\}.$$ We fix $T\in (0,T_b)$ and let $u$ be the mild solution on $[0,T]$, we define $u(t)=u(T)$ for $t\geq T$.

Now, we first pick approximating sequence $\varrho^{(n)}\in L^{p}(\mathbb{R}^{n})\cap L^{2}(\mathbb{R}^{n})$ such that $\varrho^{(n)}\geq 0$ and $\varrho_{0}^{(n)}\to \varrho_{0}$ in $L^{p}(\mathbb{R}^{n})$. For example, we can choose $f_n\in C_{c}^{\infty}(\mathbb{R}^{n})$. Denote $f_n\vee 0:=\max(f_n,0)$ and picking $\varrho_{0}^{n}=f_n\vee 0$ suffices because $\vert u_0-f_n\vee 0 \vert \leq \vert u_0-f_n \vert$ due to the fact $u_0\geq 0.$

We consider the following problem
\begin{equation}\label{111}
    \varrho^{(n)}=S_{\alpha}^{\beta}u_0-\int_{0}^{t}T_{\alpha}^{\beta}(t-s)\left(\nabla_{x}((\varrho^{(n)})^{+}\nabla v)+a\varrho^{(n)}-b{\varrho^{(n)}}^{2}\right)(s)ds.
\end{equation}
We denote
$$a_{n}(x,t)=u_{n}*J_{\varepsilon},$$
(hence $\Delta v_{n}$) is a smooth function with derivatives bounded in $[0,\infty)\times \mathbb{R}^{n}$.

For $t>0$, since $E_{\beta}(-s)\sim C_{1}s^{-1}$ as $s\to \infty$, we have
\begin{equation}
    \begin{split}
        \Vert S_{\alpha}^{\beta} u \Vert_{H^{\sigma,2}}^{2}&=\int_{\mathbb{R}^{n}}(1+\lvert \xi \rvert^{2\sigma})(E_{\beta}(-t^{\beta}\lvert \xi \rvert^{\alpha}))^{2}\lvert \hat{u}_{\varepsilon}\rvert^{2}d\xi\\
        &\leq C\Vert u \Vert_{L^2}^{2}t^{-\frac{2\sigma\beta}{\alpha}}, \ \sigma\in [0,\alpha]
    \end{split}
\end{equation}
To deal with the second term in (\ref{111}), due to the definition of $T_{\alpha}^{\beta}$, we note that for any norm, the following inequality holds
\begin{equation}
    \left \Vert \int_{0}^{t}T_{\alpha}^{\beta}\nabla \cdot Bds \right\Vert\leq \int_{0}^{t}(t-s)^{\alpha-1}\Vert E_{\beta,\beta}(-(t-s)^{\beta}A)\nabla \cdot B(s)\Vert ds.
\end{equation}
Hence, aiming to compute $H^{\sigma,2}$ norm, from \cite{li2018cauchy}, we obtain
\begin{equation}
    \Vert E_{\beta,\beta}(-(t-s)^{\beta}A)\nabla\cdot B(s)\Vert_{H^{\sigma,2}}^{2} \leq C(1+(t-s)^{-\frac{\beta(2\sigma+2-2\delta)}{\alpha}})\Vert B\Vert_{H^{\delta,2}}^{2}.
\end{equation}
Since $E_{\beta,\beta}(-s)\to CS^{-2}$ as $s\to \infty$, the above inequality is valid if $2\sigma-2+2\delta\leq 4\alpha$. This computation implies that
\begin{equation}
    \left\Vert \int_{0}^{T}(t-s)\nabla\cdot Bds\right\Vert_{H^{\sigma,2}}\leq C\left((t-s)^{\beta-1-\beta\frac{2\sigma+2-2\delta}{2\alpha}}\right)\Vert B\Vert_{H^{\delta,2}}.
\end{equation}
We find that
\begin{equation}
    \varrho^{(n)}\in C^{\beta}(0,\infty;L^{2}(\mathbb{R}^{n}))\cap C^{\infty}(0,\infty;L^{2}(\mathbb{R}^{n})).
\end{equation}
Using the time regularity, we have that in $C(0,\infty; H^{-\alpha,2}(\mathbb{R}^{n}))$, $\varrho^{(n)}$ solves the following problem in strong sense:
\begin{equation}
\left\{
\begin{aligned}
    & \partial_{t}^{\beta}\varrho^{(n)}+ (-\Delta)^{-\frac{\alpha}{2}}\varrho^{(n)} = - \nabla \cdot ((\varrho^{(n)})^{+}a_{n}(x,t)) + a\varrho^{(n)}-b(\varrho^{(n)})^{2},  \\
    & \varrho(x,0)=\varrho_{0}^{(n)}(x)\geq 0
\end{aligned}
\right.
\end{equation}
with the time regularity and Lemma \ref{2.8.1}, we find that in $H^{-\alpha,2}(\mathbb{R}^{n})$
\begin{equation}\label{111.1}
    \partial_{t}^{\beta}\varrho^{(n)}=\frac{1}{\Gamma(1-\beta)}\left(\frac{\varrho^{(n)}-\varrho_{0}^{(n)}}{t^{\beta}}+\beta\int_{0}^{t}\frac{\varrho^{(n)}-\varrho^{(n)}(x,s)}{(t-s)^{\beta+1}}ds\right).
\end{equation}
It is clear that
$$u\mapsto a_{n}u^{+}+au-bu^{2}$$
is bounded in $L^{2}(\mathbb{R}^{n})$ and $H^{1,2}(\mathbb{R}^{n})$. By interpolation, this mapping id bounded in $H^{\delta,2}$, $0\leq \delta \leq 1.$

We now denote
$$V_{1}^{n}=(\varrho^{(n)})^{+}a_{n},\ V_{2}^{n}=a\varrho^{(n)},\ V_{3}^{n}=b{\varrho^{(n)}}^{2}.$$
We know $u_{n}\in C([0,T];L^{2}(\mathbb{R}^{n}))$, hence we can pick $\delta=0$, $2\sigma+2+2\delta<2\alpha$ or $\sigma<2\alpha-1$, we find that
\begin{equation}
\begin{split}
    \Vert u_{n}\Vert_{H^{\sigma,2}}&\leq Ct^{-\sigma\beta}{\alpha}+C(T)\left(1+\int_{0}^{t}(t-s)^{\beta-1-\beta\frac{2\sigma+2-2\delta}{2\alpha}}\right)\Vert V_{1}^{n}+V_{2}^{n}+V_{3}^{n}\Vert_{L^{2}}\\
    &\leq Ct^{-\frac{\sigma\beta}{\alpha}}.
    \end{split}
\end{equation}
Hence, for $t>0$
$$\varrho^{(n)}\in H^{\sigma,2}(\mathbb{R}^{n}), \ \sigma\in [0,\alpha).$$
Consequently, for $t>0$, $\varrho^{(n)}\in H^{1,2}(\mathbb{R}^{n})$ and $(-\Delta)^{\frac{\alpha}{2}}\in H^{-\varepsilon,2}(\mathbb{R}^{n})$ for any $\varepsilon >0.$

The right hand side of (\ref{111.1}) also makes sense in $C(0,\infty;L^{2}(\mathbb{R}^{n}))$. Since $(\varrho^{(n)})^{-}\in L^{2}(\mathbb{R}^{n})$, we can multiply $(\varrho^{(n)})^{-}$ on both sides of (\ref{111.1}) for $t>0$ and take integral with respect to $x$. Together with the facts $\varrho^{(n)}\in H^{1,2}(\mathbb{R}^{n})$ and $(-\Delta)^{\frac{\alpha}{2}}\varrho^{(n)}\in H^{1,2}(\mathbb{R}^{n})$ for $t>0$, we have
\begin{equation}
    \begin{split}
        \int_{\mathbb{R}^{n}}\partial_{t}^{\beta}\varrho^{(n)}(\varrho^{(n)})^{-}dx+\langle (-\Delta)^{\frac{\alpha}{2}}\varrho^{(n)}, (\varrho^{(n)})^{-} \rangle &=-\int_{\mathbb{R}^{n}}\nabla\cdot ((\varrho^{(n)})^{+}a_{n}(x,t)(\varrho^{(n)})^{-}dx\\
        &+\int_{\mathbb{R}^{n}}\langle a\varrho^{(n)}, (\varrho^{(n)})^{-}\rangle dx\\
        &-\int_{\mathbb{R}^{n}}\langle b{\varrho^{(n)}}^{2}, (\varrho^{(n)})^{-}\rangle dx.
    \end{split}
\end{equation}
(i) 
$$\langle (-\Delta)^{\frac{\alpha}{2}}\varrho^{(n)}, (\varrho^{(n)})^{-}\rangle\leq 0,$$
(ii) 
$$-\langle \nabla \cdot \left((\varrho^{(n)})^{+}a_{n}(x,t)\right) ,(\varrho^{(n)})^{-}\rangle=0.$$
The proof process for (i) and (ii) is similar to the method in \cite{li2018cauchy}. 

The same can be obtained
$$\langle a\varrho^{(n)}, (\varrho^{(n)})^{-}\rangle \leq 0,$$
$$\langle b{\varrho^{(n)}}^{2}, (\varrho^{(n)})^{-}\rangle \leq 0.$$
From the above, it can be obtained that for $t>0$
\begin{equation}
    0\leq -\frac{1}{2}\left(\partial_{t}^{\beta}\Vert (\varrho^{(n)})^{-}\Vert_{L^{2}}^{2}\right).
\end{equation}
Since $\Vert ((\varrho^{(n)}))^{-}\Vert_{L^{2}}^{2}$ is continuous in time. $(\varrho^{n})^{-}=0$ follows from Lemma \ref{2.8.2}. This means that $\varrho^{(n)}$ indeed satisfies the following equation
\begin{equation}
    \varrho^{(n)}(t)=S_{\alpha}^{\beta}\varrho_{0}^{(n)}-\int_{0}^{T}T_{\alpha}^{\beta}(t-s)\left(\nabla \cdot (\varrho^{(n)})\nabla v_{n}+a\varrho^{(n)}-b{\varrho^{(n)}}^{2}\right)ds
\end{equation}
and $\varrho^{(n)}\geq 0$.

Now that $\varrho_{0}^{(n)}\in L^{p}(\mathbb{R}^{n})$, since $J_{\varepsilon}(x)*u_{n}$ is smooth and bounded, we find that $\varrho^{(n)}\in C(0,\infty;L^{p}(\mathbb{R}^{n}))$. Then for $t\in [0,T]$,
\begin{equation}
    \begin{split}
        \Vert \varrho^{(n)}-u(t)\Vert_{L^{p}}&\leq \Vert \varrho_{0}^{(n)}-u(t)\Vert_{L^{p}}\\
        &+\int_{0}^{t}(t-s)^{-\frac{n\beta}{\alpha}\left(\frac{1}{q}-\frac{1}{p}\right)-\frac{\beta}{\alpha}+\beta-1}\left(
        {\begin{array}{l}
        \Vert \varrho^{(n)}\nabla v_{n}-u\nabla v\Vert_{L^{q}}\\
        +a\Vert \varrho^{(n)}-u\Vert_{L^{q}}\\
        +b\Vert {\varrho^{(n)}}^{2}-u^{2}\Vert_{L^{q}}\\
         \end{array} }
        \right)ds\\
        &\leq \Vert \varrho_{0}^{(n)}-u(t)\Vert_{L^{p}}\\
        &+\int_{0}^{t}(t-s)^{-\frac{n\beta}{\alpha}\left(\frac{1}{q}-\frac{1}{p}\right)-\frac{\beta}{\alpha}+\beta-1}\Vert \varrho^{(n)}-u\Vert_{L^{p}}\Vert \nabla v_{n}\Vert_{L^{\frac{pq}{p-q}}}ds\\
        &+\int_{0}^{t}(t-s)^{-\frac{n\beta}{\alpha}\left(\frac{1}{q}-\frac{1}{p}\right)-\frac{\beta}{\alpha}+\beta-1}\Vert u \Vert_{L^{p}}\Vert \nabla v_{n}-\nabla v\Vert_{L^{\frac{pq}{p-q}}}ds\\
        &+a\int_{0}^{t}(t-s)^{-\frac{n\beta}{\alpha}\left(\frac{1}{q}-\frac{1}{p}\right)-\frac{\beta}{\alpha}+\beta-1}\Vert \varrho^{(n)}-u\Vert_{L^{q}}ds\\
        &+b\int_{0}^{t}(t-s)^{-\frac{n\beta}{\alpha}\left(\frac{1}{q}-\frac{1}{p}\right)-\frac{\beta}{\alpha}+\beta-1}\Vert {\varrho^{(n)}}^{2}-u^{2}\Vert_{L^{q}}ds.
    \end{split}
\end{equation}
From (\ref{11111}), we can similarly obtain 
\begin{equation}
    \Vert \varrho^{(n)}(t)-u(t)\Vert_{L^{p}}\leq C\int_{0}^{t}(t-s)^{-\frac{n\beta}{\alpha}\left(\frac{1}{q}-\frac{1}{p}\right)-\frac{\beta}{\alpha}+\beta-1}\Vert \varrho^{(n)}-u\Vert_{L^{p}}ds+\delta_{n}(t),
\end{equation}
where $\delta_{n}(T)=\Vert \varrho_{0}^{(n)}-u_{0}\Vert_{L^{p}}\to 0$, as $n\to\infty$.

By the comparison principle (Lemma \ref{2.8.3}), we have 
\begin{equation}
    \Vert \varrho^{(n)}(t)-u(t)\Vert_{L^{p}}\leq \rho(t),
\end{equation}
where $\rho(t)$ solves the equation
\begin{equation}
\partial_{t}^{-\frac{n\beta}{\alpha}\left(\frac{1}{q}-\frac{1}{p}\right)-\frac{\beta}{\alpha}+\beta-1}\rho_{n}(t)=C\rho_{n}(t),\ u_{n}(0)=\delta_{n}(T).
\end{equation}
It is easy to see that $\rho_{n}(t)\to 0$, as $n\to \infty$ for $t\in[0,T]$, we have
\begin{equation}
    \varrho^{(n)}(t)\to u(t),\ \ {\rm in\ }C(0,T;L^{p}(\mathbb{R}^{n})).
\end{equation}
Since $\varrho^{(n)}(t)\geq 0$, we then have $u(t)\geq 0$ for $t\in [0,T]$. Since $T\in (0,T_{b})$is arbitrary, the claim follows.\qed
\end{pro}
\begin{remark}
Different from literature \cite{jiang2021weak}, we consider the nonlinear convection-diffusion equation and prove the non-negativity of its mild solution.
\end{remark}

Before proving the global existence of weak solutions, we start with the construction of some a prior estimates of the solutions to the mollified Equ. (\ref{2.13}).
\begin{lemma}\label{1.1.4}
For $\gamma\in \mathbb{R}$ and $\beta \in(0,1)$. Suppose $h(x,t)$ is smooth and all the derivatives of $h(x,t)$ and bounded. Then the following holds.
\begin{enumerate}[(i)]
        \item If $u_{0}\in L^{1}(\mathbb{R}^{n})\cap H^{\gamma}(\mathbb{R}^{n})$, then for any $T>0$, (\ref{2.13}) has a unique mild solution in $C\left([0,T];H^{\gamma}(\mathbb{R}^{n})\right)$, such that the nonnegativity in Lemma \ref{3.1.1.1.1.1} hold true.
        \item For the unique mild solution in (i) and any $T>0$,
        $$u\in C^{0,\beta}\left(0,T;H^{\gamma}(\mathbb{R}^{n})\right)\cap C^{\infty}\left(0,T;H^{\gamma}(\mathbb{R}^{n})\right).$$
        Moreover, the following holds strongly in $C\left([0,T;H^{\gamma-\alpha}(\mathbb{R}^{n})\right)$,
        $$\partial_{t}^{\alpha}u=\frac{1}{\Gamma(\alpha)}\int_{0}^{t}(t-s)^{\beta-1}\frac{\partial u(s)}{\partial s}ds=-\nabla \cdot (uh(x,t))-(-\Delta)^{\frac{\alpha}{2}}+F.$$
    \end{enumerate}
\end{lemma}
\begin{remark}
The difference from the linear convective-diffusion equation in \cite{jiang2021weak} is that we prove the existence of mild solutions to the nonlinear convective-diffusion equation. The difficulty with this proof lies in our treatment of nonlinear terms.
\end{remark}
\begin{pro}
We take $F(x,t)=au(x,t)-bu^{2}(x,t)$. Since $E'_{\beta}(z)$ is an analytic function in the whole complex plane $z\in \mathbb{C}$ and
$$E'_{\beta}(-s)\sim -C_{0}s^{-2}, \ \ {\rm as\ } s\to \infty,$$ 
we conclude that 
$$\underset{0\leq s\leq \infty}{\sup}E'_{\beta}(-s)s^{\sigma}\leq C,\ \forall \sigma \leq 2.$$
Consequently, Plancherel theorem  imply that
\begin{equation}
    \begin{split}
        \Vert E'_{\beta}(-t^{\beta}A)\nabla f\Vert_{H^{\gamma}(\mathbb{R}^{n})}^{2}&=\Vert \mathcal{F}^{-1}\left((1+\lvert \xi\rvert^{2})^{\frac{\gamma}{2}}\mathcal{F}(E'_{\beta}(-t^{\beta}A)\nabla f)\right) \Vert_{2}^{2}\\
        &=\Vert (1+\lvert \xi\rvert^{2})^{\frac{\gamma}{2}}(E'_{\beta}(-t^{\beta}\lvert \xi \rvert^{\alpha})\lvert \xi \rvert \hat{f}_{\xi})\Vert_{2}^{2}\\
        &\leq C\tau^{-\frac{2\beta}{\alpha}}\int_{\mathbb{R}^{n}}\lvert \hat{f}_{\xi} \rvert^{2}(1+\lvert \xi\rvert^{2})^{\gamma}ds\\
        &=Ct^{-\frac{2\beta}{\alpha}}\Vert f\Vert_{H^{\gamma}(\mathbb{R}^{n})}^{2}.
    \end{split}
\end{equation}
We construct the iterative sequence
\begin{equation}\label{3.1.1.1}
    u_{0}(t)=u_{0}
\end{equation}
\begin{equation}\label{3.1.1.2}
    \begin{split}
        u_{n}(t)&=E_{\beta}(-t^{\beta}A)u_{0}\\
        &-\beta\int_{0}^{t}\tau^{\beta-1}E'_{\beta}(-\tau^{\beta}A)\left(-\nabla \cdot (u_{n-1}h)|_{t-\tau}+(au_{n-1}-bu_{n-1}^{2})|_{t-\tau}\right)d\tau.
    \end{split}
\end{equation}
Fix $T>0$ and define $E_{n}(t)=u_{n}(t)-u_{n-1}(t)$. Combine (\ref{3.1.1.1}) and (\ref{3.1.1.2}), compute directly to yield
\begin{equation}
    \begin{split}
        u_{1}&=E_{\beta}(-t^{\beta}A)u_{0}\\
        &-\beta\int_{0}^{t}\tau^{\beta-1}E'_{\beta}(-\tau^{\beta}A)\left(-\nabla \cdot (u_{1}h)|_{t-\tau}+(au_{1}-bu_{1}^{2})|_{t-\tau}\right)d\tau\\
        &\leq E_{\beta}(-t^{\beta}A)u_{0}-\beta\int_{0}^{t}\tau^{\beta-1}E'_{\beta}(-\tau^{\beta}A)\left(-\nabla \cdot (u_{1}h+au_{1})|_{t-\tau}\right)d\tau.
    \end{split}
\end{equation}
It is easy to see that
\begin{equation}\label{3.1.1.3}
    \begin{split}
        \Vert u_1 \Vert_{C(0,T;H^{\gamma}(\mathbb{R}^{n}))}&\leq \Vert u_{0}\Vert_{C(0,T;H^{\gamma}(\mathbb{R}^{n}))}\left(1+C_{1}\beta\int_{0}^{t}\tau^{\beta-\frac{\beta}{\alpha}-1}(1+\tau^{\frac{\beta}{\alpha}})d\tau\right)\\
        &\leq \Vert u_{0}\Vert_{C(0,T;H^{\gamma}(\mathbb{R}^{n}))}\left(1+\frac{\alpha}{\alpha-1}C_1T^{\beta-\frac{\beta}{\alpha}}+C_1T^{\beta}\right).
    \end{split}
\end{equation}
Therefore, for any $0\leq t \leq T$, (\ref{3.1.1.3}) shows that 
\begin{equation}
    \Vert E_{1} \Vert_{C(0,T;H^{\gamma}(\mathbb{R}^{n}))}\leq M,
\end{equation}
where $M=(2+\frac{\alpha}{\alpha+1})T^{\beta-\frac{\beta}{\alpha}}$.

By Lemma \ref{3.1.1.1.1.1}, the induction formula reads
\begin{equation}
    \begin{split}
        E_{n}(t)&=\beta\int_{0}^{t}\tau^{\beta-1}E'_{\beta}(-\tau^{\beta}A)\left(-\nabla \cdot (E_{n-1}h+aE_{n-1}-b(u_{n}^{2}-u_{n-1}^{2}))|_{t-\tau}\right)d\tau\\
        &\leq C_2\beta\int_{0}^{t}\tau^{\beta-1}E'_{\beta}(-\tau^{\beta}A)\left(-\nabla \cdot (E_{n-1}h+aE_{n-1}-b(u_{n}^{2}-u_{n-1}^{2}))|_{t-\tau}\right)d\tau\\
        &\leq C_2\beta\int_{0}^{t}\tau^{\beta-1}E'_{\beta}(-\tau^{\beta}A)\left(-\nabla \cdot (E_{n-1}h+aE_{n-1})|_{t-\tau}\right)d\tau.
    \end{split}
\end{equation}
As a result,
\begin{equation}
    \begin{split}
        \Vert E_{n}(t) \Vert_{C(0,T;H^{\gamma}(\mathbb{R}^{n}))}&\leq C_2\beta\underset{0\leq z \leq t}{\sup}\int_{0}^{z}\left(\tau^{\beta-\frac{\beta}{\alpha}-1}+\tau^{\beta-1}\right)\Vert E_{n-1}(t) \Vert_{C(0,T;H^{\gamma}(\mathbb{R}^{n}))}d\tau\\
        &= C_2\beta\int_{0}^{t}\left(\tau^{\beta-\frac{\beta}{\alpha}-1}+\tau^{\beta-1}\right)\Vert E_{n-1}(t) \Vert_{C(0,T;H^{\gamma}(\mathbb{R}^{n}))}d\tau\\
        &\leq C\beta g_{\beta-\frac{\beta}{\alpha}}\ast \Vert E_{n-1}(t) \Vert_{C(0,T;H^{\gamma}(\mathbb{R}^{n}))}. 
    \end{split}
\end{equation}
From the above formula, one has
\begin{equation}\label{3.1.1.4}
    \Vert E_{2}(t) \Vert_{C(0,T;H^{\gamma}(\mathbb{R}^{n}))}\leq CMg_{\beta-\frac{\beta}{\alpha}+1}(t).
\end{equation}
and
\begin{equation}\label{3.1.1.5}
     \Vert E_{n}(t) \Vert_{C(0,T;H^{\gamma}(\mathbb{R}^{n}))}\leq C^{n-1}Mg_{(n-1)(\beta-\frac{\beta}{\alpha})+1}(t).
\end{equation}
Estimates (\ref{3.1.1.4}) and (\ref{3.1.1.5}) indicate that $u=u_0+\sum_{n=1}^{\infty}E_{n-1}$ converges in $C([0,T];H^{\gamma}(\mathbb{R}^{n}))$. In other words, $u^{n}\to u$ in $C(0,T;H^{\gamma}(\mathbb{R}^{n}).$ Hence, $u$ is a mild solution.

The uniqueness follows in a similar way. Consider two mild solutions $u_1$ and $u_2$. Define $\varsigma=u_1-u_2$. Then,
\begin{equation}
    \begin{split}
        \varsigma&=\beta\int_{0}^{t}\tau^{\beta-1}E'_{\beta}(-\tau^{\beta}A)\left(-\nabla \cdot (\varsigma h+a\varsigma-b(u_{1}^{2}-u_{2}^{2}))|_{s}\right)d\tau\\
        &\leq C\beta\int_{0}^{t}\tau^{\beta-1}E'_{\beta}(-\tau^{\beta}A)\left(-\nabla \cdot (\varsigma h+a\varsigma-b\varsigma)|_{s}\right)d\tau.
    \end{split}
\end{equation}
We have
\begin{equation}
    \Vert \varsigma \Vert_{C(0,T;H^{\gamma}(\mathbb{R}^{n}))}\leq C\underset{0\leq z \leq t}{\sup}\int_{0}^{z}\left(\tau^{\beta-\frac{\beta}{\alpha}-1}+\tau^{\beta-1}\right)\Vert \varsigma \Vert_{C(0,T;H^{\gamma}(\mathbb{R}^{n}))}d\tau.
\end{equation}
The above formula and the comparison principle in 9 Proposition 5 \cite{feng2018continuous} ] imply that $\Vert \varsigma \Vert_{C(0,T;H^{\gamma}(\mathbb{R}^{n}))}=0$ for $0\leq t\leq T$. Thus, the uniqueness is establishes.

(ii) Assume that $u(x,t)$ is the mild solution which satisfies
\begin{equation}
    \begin{split}
        u(x,t)&=E_{\beta}(-t^{\beta}A)u_{0}\\
        &-\beta\int_{0}^{t}(t-s)^{\beta-\frac{\beta}{\alpha}-1}\left((t-s)^{\frac{\beta}{\alpha}}E'_{\beta}(-\tau^{\beta}A)\left(-\nabla \cdot (uh+au-bu^2)|_{t-\tau}\right)\right)d\tau.
    \end{split}
\end{equation}
Notice that $(t-s)^{\frac{\beta}{\alpha}}E'_{\beta}(-t^{\beta}A)\nabla: H^{\gamma}(\mathbb{R}^{n})\to  H^{\gamma}(\mathbb{R}^{n})$ is a bounded operator, $E_{\beta}(-t^{\beta}A)u_0\in C^{\beta}(0,T;H^{\gamma}(\mathbb{R}^{n}))$ is proved in \cite{taylorremarks}, according to \cite{li2018some}, for $T>0$,
\begin{equation}
    u\in C^{0,\beta}([0,T];H^{\gamma}(\mathbb{R}^{n}))\cap C^{\infty} ([0,T];H^{\gamma}(\mathbb{R}^{n})).
\end{equation}
Since $E_{\beta}(-t^{\beta}A)\varphi$ solves the fractional diffusion equation $\partial_{t}^{\alpha}u+(-\Delta)^{\frac{\alpha}{2}}u=0$, we deduce
\begin{equation}\label{3.1.1.7}
    E_{\beta}(-t^{\beta}A)\varphi=\varphi-\frac{1}{\Gamma(\beta)}\int_{0}^{t}(t-s)^{\beta-1}AE_{\beta}(-s^{\beta}A)\varphi ds,
\end{equation}
taking the derivative on $t$, we find the operator identity,
\begin{equation}\label{3.1.1.8}
\begin{split}
    -\beta t^{\beta-1}E'_{\beta}(-t^{\beta}A)&=-\frac{1}{\Gamma(\beta)}t^{\beta-1}I\\
    &+\frac{1}{\Gamma{\beta}}\int_{0}^{t}(t-s)^{\beta-1}As^{\beta}E'_{\beta}(-s^{\beta}A)ds.
    \end{split}
\end{equation}
From (\ref{3.1.1.7}), (\ref{3.1.1.8}) and $Au\in C(0,\infty;H^{\gamma}(\mathbb{R}^{n})$, we find that 
\begin{equation}\label{3.1.1.9}
   u\in C^{\beta}(0,T;H^{\gamma-\alpha}(\mathbb{R}^{n}))\cap C^{\infty}(0,\infty;H^{\gamma-\alpha}(\mathbb{R}^{n})),
\end{equation}
and
\begin{equation}\label{3.1.1.10}
    u(x,t)=u_0+\frac{1}{\Gamma{\beta}}\int_{0}^{t}(t-s)^{\beta-1}\left(-Au(s)-\nabla \cdot (uh)(s)\right)ds.
\end{equation}
Due to the time regularity estimates (\ref{3.1.1.9}) and (\ref{3.1.1.10})
\begin{equation}
\begin{split}
    \partial_{t}^{\beta}u&=\frac{1}{\Gamma(1-\beta)}\int_{0}^{t}(t-s)^{\beta-1}\frac{\partial u(s)}{\partial s}ds\\
    &= -\nabla \cdot (uh)-(-\Delta)^{\frac{\alpha}{2}}u+au-bu^2
\end{split}
\end{equation}
to be held in $C(0,T;H^{\gamma-\alpha}(\mathbb{R}^{n})).$\qed
\end{pro}
\begin{remark}
We adopt a method similar to \cite{jiang2021weak}, combined with Lemma \ref{3.1.1.1.1.1} to prove the existence and uniqueness of the mild solution of the nonlinear convection-diffusion equation, which plays an important role in proving existence of the strong solution of the Mollified system in Theorem \ref{1.1.1} and Theorem \ref{1.1.2}. 
\end{remark}

\subsection{The case of \ \texorpdfstring {$b\geq 1$}{lg}}
\begin{lemma}\label{3.1}
Under the assumption of Theorem \ref{1.1.1}(i), the solution $(u,v)$ of system (\ref{1.1}) satisfies the following estimation
\begin{flalign}\label{3.1.1}
        && \Vert u \Vert_{L^{\infty}([0,T];L^{q}(\mathbb{R}^{n}))} &\leq C, \ \ \Vert D^{\frac{\alpha}{2}}u^{\frac{q}{2}} \Vert_{L^{2}(0,T;L^{2}(\mathbb{R}^{n}))}\leq C.\\
       && \Vert \nabla v \Vert_{L^{\infty}(0,\infty;L^{s_1}(\mathbb{R}^{n}))} &\leq C, \ \ s_1 > \frac{n}{n-1}
    \end{flalign}
\end{lemma}
\begin{pro}
It is obtained by multiplying (\ref{1.1}) with $qu^{q-1}$ and by Lemma \ref{2.9}, then integrate over $\mathbb{R}^{n}$
\begin{equation}\label{3.1.2}
\begin{split}
    \int_{\mathbb{R}^{n}}qu^{q-1}\partial_{t}^{\beta}udx
    &=-\int_{\mathbb{R}^{n}} (-\Delta)^{\frac{\alpha}{2}}qu^{q-1}dx-\int_{\mathbb{R}^{n}}\nabla \cdot (u\nabla v)qu^{q-1}dx\\
    &+ au\int_{\mathbb{R}^{n}}u^{q}dx-bu\int_{\mathbb{R}^{n}}u^{q+1}dx\\
    &\leq -\frac{4(q-1)}{q}\int_{\mathbb{R}^{n}}\lvert D^{\frac{\alpha}{2}}u^{\frac{q}{2}}\rvert^{2}dx-(q-1)\int_{\mathbb{R}^{n}}u^{q}\Delta v dx\\
    &+ aq\int_{\mathbb{R}^{n}}u^{q}dx-bq\int_{\mathbb{R}^{n}}u^{q+1}dx.
    \end{split}
\end{equation}
Substitute $-\Delta v = u$ into (\ref{3.1.2}) to get
\begin{equation}\label{3.1.3}
\begin{split}
    q\int_{\mathbb{R}^{n}}u^{q-1}\partial_{t}^{\beta}udx+\frac{4(q-1)}{q}\int_{\mathbb{R}^{n}}\lvert D^{\frac{\alpha}{2}}u^{\frac{q}{2}}\rvert^{2}dx&\leq \int_{\mathbb{R}^{n}}(q-1-bq)u^{q+1}dx+aq\int_{\mathbb{R}^{n}}u^{q}dx.
    \end{split}
\end{equation}
$q-1-bq<0$ from $b\geq 1$ and by Lemma \ref{2.10}, we obtain
\begin{equation}\label{3.1.4}
\partial_{t}^{\beta}\int_{\mathbb{R}^{n}}u^{q}dx+\frac{4(q-1)}{q}\int_{\mathbb{R}^{n}}\lvert D^{\frac{\alpha}{2}}u^{\frac{q}{2}}\rvert^{2}dx\leq aq\int_{\mathbb{R}^{n}}u^{q}dx.
\end{equation}
It is easy to see that $\partial_{t}^{\beta}\int_{\mathbb{R}^{n}}u^{q}dx\leq aq\int_{\mathbb{R}^{n}}u^{q}dx.$ By Lemma \ref{2.11}, we have
$$\int_{\mathbb{R}^{n}}u^{q}dx\leq u_{0}E_{\alpha}(aqt^{\alpha})\leq C.$$
Thai is 
$$\Vert u \Vert_{L^{\infty}(0,T;L^{q}(\mathbb{R}^{n}))}\leq C.$$
It is easy to know from (\ref{3.1.4})
$$\frac{4(q-1)}{q}\int_{\mathbb{R}^{n}}\lvert D^{\frac{\alpha}{2}}u^{\frac{q}{2}}\rvert^{2}dx\leq C.$$
So we have $D^{\frac{\alpha}{2}}u^{\frac{q}{2}}\in L^{2}(0,T;L^{2}(\mathbb{R}^{n}))$, that is
$$\Vert D^{\frac{\alpha}{2}}u^{\frac{q}{2}} \Vert_{L^{2}(0,T;L^{2}(\mathbb{R}^{n}))}\leq C.$$
Consequently, in view of (\ref{3.1.4}), we obtain that $u\in L^{\infty}(0,T;L^{q}(\mathbb{R}^{n}))\cap L^{2}(0,T;H^{\frac{\alpha}{2}}(\mathbb{R}^{n}))$.

We know that
\begin{equation}
    \begin{split}
        \Vert \nabla v \Vert_{L^{s_1}} & =\Vert F* u\Vert\\
        &=-\left\Vert \frac{C_{*}x}{\lvert x\rvert^{n}}* u\right\Vert_{L^{s_1}}\\
        &\leq C\left\Vert \frac{x}{\lvert x\rvert^{n}}* u \right\Vert_{L^{s_1}}\\
        &\leq C\Vert \lvert x \rvert^{-(n-1)}\Vert_{L_{\omega}^{p}}\Vert u \Vert_{L^{q}}\\
        &\leq C\Vert u \Vert_{L^{q}},
    \end{split}
\end{equation}
where $q>1$, $p=\frac{n}{n-1}$, and satisfies $1+\frac{1}{s_1}=\frac{1}{p}+\frac{1}{q}$, i.e., $s_{1}=\frac{nq}{n-q}>\frac{n}{n-1}$.
\end{pro}

\begin{lemma}\label{3.2.1.1.2}
Under the assumption of Theorem \ref{1.1.1}(i), the time derivative of $u$ satisfies the following estimation
\begin{equation}
    \Vert \partial_{t}^{\beta}u\Vert_{L^{q_{1}}(0,T;W^{-\alpha,\frac{nq}{2n-q}})}\leq C,\ q_{1}>1
\end{equation}
\end{lemma}
\begin{pro}
Since $u$ is uniformly bounded in $L^{1}(\mathbb{R}^{n})\cap L^{\infty}(\mathbb{R}^{n})$, it is also in $L^{\frac{nq}{2n-q}}(\mathbb{R}^{n})$ for $1\leq \frac{nq}{2n-q}\leq q <\infty$. For any fixed $T>0$, we take the test function $\varphi(x,t)$ with $p_{1}>1$, 
\begin{equation}\label{3.2.1.1}
\Vert \varphi \Vert_{L^{p_{1}}\left((0,T);W^{\alpha,\frac{nq}{nq-2n+q}}(\mathbb{R}^{n})\right)}\leq 1.
\end{equation}
We have
\begin{equation}
    \begin{split}
        \int_{\mathbb{R}^{n}}\partial_{t}^{\beta}u\varphi dx&=-\int_{\mathbb{R}^{n}}(-\Delta)^{\alpha}u\varphi dx+\int_{\mathbb{R}^{n}}u\nabla v \cdot\nabla \varphi dx\\
        &+a\int_{\mathbb{R}^{n}}u\varphi dx-b\int_{\mathbb{R}^{n}}u^{2}\varphi dx.
    \end{split}
\end{equation}
Next we estimate $\langle \partial_{t}^{\beta}u, \varphi \rangle$
\begin{equation}
\begin{split}
    \lvert \langle \partial_{t}^{\beta}u, \varphi \rangle \rvert 
    & \leq \lvert - \langle u, D^{\alpha}\varphi \rangle \rvert +\lvert \langle u\nabla v, \nabla \varphi \rangle \rvert + \lvert \langle au,  \varphi \rangle \rvert+ \lvert \langle bu^{2},  \varphi \rangle \rvert\\
    &\leq \Vert u \Vert_{L^{\infty}((0,T);L^{q}(\mathbb{R}^{n}))} \Vert \nabla v \Vert_{L^{\infty}((0,T);L^{\frac{nq}{2n-q}}(\mathbb{R}^{n}))}\int_{0}^{T}\Vert \nabla \varphi \Vert_{L^{\frac{nq}{q-2n+q}}(\mathbb{R}^{n})}dt\\
    &+ \Vert u\Vert_{L^{\infty}\left((0,T);L^{\frac{nq}{2n-q}}(\mathbb{R}^{n})\right)}\int_{0}^{T}\Vert D^{\alpha}\varphi\Vert_{L^{\frac{nq}{nq-2n+q}}}dt\\
    &+a\Vert u \Vert_{L^{\infty}\left((0,T);L^{\frac{nq}{2n-q}}(\mathbb{R}^{n})\right)}\Vert \varphi\Vert_{L^{\frac{nq}{nq-2n+q}}}dt\\
    &+b\Vert u^{2} \Vert_{L^{\infty}\left((0,T);L^{\frac{nq}{2n-q}}(\mathbb{R}^{n})\right)}\Vert \varphi\Vert_{L^{\frac{nq}{nq-2n+q}}}dt.
\end{split}
\end{equation}
According to (\ref{3.2.1.1}), for any $T>0$
\begin{align}
    & \int_{0}^{T}\Vert \partial_{t}^{\beta}u\Vert_{L^{q_{1}}\left((0,T);W^{-\alpha, \frac{nq}{2n-q}}(\mathbb{R}^{n})\right)}dt\notag\\
    \leq{} &  C
    \left(
    \begin{array}{cc}
    \int_{0}^{T} \Vert u \Vert_{L^{\infty}\left((0,T);L^{\frac{nq}{2n-q}}(\mathbb{R}^{n})\right)}dt
    &\\
    +  \int_{0}^{T}\Vert u \Vert_{L^{\infty}((0,T);L^{q}(\mathbb{R}^{n}))}\Vert \nabla v \Vert_{L^{\infty}((0,T);L^{\frac{nq}{2n-q}}(\mathbb{R}^{n}))}dt   &  \\
      +\int_{0}^{T}\Vert u\Vert_{L^{\infty}\left((0,T);L^{\frac{nq}{2n-q}}(\mathbb{R}^{n})\right)}dt
      &\\
      +\int_{0}^{T}\Vert u\Vert_{L^{\infty}\left((0,T);L^{\frac{2nq}{2n-q}}(\mathbb{R}^{n})\right)}^{2}   & 
    \end{array}
    \right).\notag\\
\end{align}
Also, we have $\nabla v=C_{*}\frac{x}{\vert x \vert^{n}}* u$, and use the Hardy-Littlewood-Sobolev inequality, we obtain
\begin{equation}
    \begin{split}
        \int_{0}^{T}\Vert \nabla v \Vert_{L^{\frac{nq}{n-q}}(\mathbb{R}^{n})}dt&\leq C\left\Vert \frac{1}{\lvert x \Vert^{n-1}}* u\right\rvert_{L^{\frac{nq}{n-q}}(\mathbb{R}^{n})}\\
        &\leq C \Vert u \Vert_{L^{q}}.
    \end{split}
\end{equation}
Therefore, it is easy to see that
\begin{equation}
    \Vert \partial_{t}^{\beta} u\Vert_{L^{q_1}\left((0,T);W^{-\alpha,\frac{nq}{2n-q}}(\mathbb{R}^{n})\right)}\leq C(q_1,q,T),
\end{equation}
where $q_1>1$ satifies $\frac{1}{q_1}+\frac{1}{p_1}=1$.\qed
\end{pro}

 With the help of Lemma \ref{1.1.4} regarding fractional advection diffusion Equation (\ref{2.13}) and generalized strong compactness criteria given in Lemma \ref{3.1.1.1.1}.\\
\noindent $\mathbf{Proof\ of\ Theorem\ref{1.1.1}(i)}$:\\
we follow the method in \cite{bian2013dynamic} by taking a cutoff function $0<\varphi_{1}(x)\leq 1$, $\varphi_{1}(x)\in C_{0}^{\infty}(\mathbb{R}^{n})$, which satisfies
\begin{equation}
    \varphi_{1}(x)=\begin{cases}
    1, &{\rm if\ } \lvert x\rvert \leq 1,\\
    0,&{\rm if\ } \lvert x\rvert \leq 2,
\end{cases}
\end{equation}
Define $\varphi_{R}(x):=\varphi_{1}(\frac{x}{R})$, then $\varphi_{R}(x)\to 1$ as $R\to \infty$, 
\begin{equation}
    \begin{split}
        D^{\alpha}\varphi_{R}(x)&=C_{\alpha,n}P.V.\int_{\mathbb{R}^{n}}\frac{\varphi_{1}(\frac{x}{R})-\varphi_{1}(\frac{y}{R})}{\lvert x-y\rvert^{n+\alpha}}dy\\
        &=C_{\alpha,n}P.V.\int_{\mathbb{R}^{n}}\frac{\varphi_{1}(\frac{x}{R})-\varphi_{1}(y)}{R^{\alpha}\lvert \frac{x}{R}-y\rvert^{n+\alpha}}dy=\frac{1}{R^{\alpha}}D^{\alpha}\varphi_{1}(\frac{x}{R}),
    \end{split}
\end{equation}
where $P.V.$ denotes the Cauchy principle value.\\
The proof can be divided into 5 steps.\\
\textbf{Step 1.} (Regularization and a priori estimates)\\
To prove the existence of weak solutions with the above properties, we consider the following regularization problem for $\varepsilon >0$
\begin{equation}\label{3.1.5}
\left\{
\begin{aligned}
    & \partial_{t}^{\beta} u_{\varepsilon} = -(-\Delta)^{\frac{\alpha}{2}} u_{\varepsilon}) - \nabla \cdot (u_{\varepsilon}\nabla v_{\varepsilon}) + au_{\varepsilon} -bu_{\varepsilon}^{2},  & x\in \mathbb{R}^{n}, t>0\\
    & -\Delta v_{\varepsilon} = J_{\varepsilon} * u_{\varepsilon}, & x\in \mathbb{R}^{n}, t>0\\
    & u_{\varepsilon}(x,0)  = u_{\varepsilon_0}, & x\in \mathbb{R}^{n},\\
\end{aligned}
\right.
\end{equation}
where $J_{\varepsilon}(x) = \frac{1}{\varepsilon^{n}}J(\frac{x}{\varepsilon}), J(x) = \frac{1}{\alpha(n)}\left(1+\lvert x \lvert^{2}\right)^{-\frac{n+2}{2}}$ satisfying $\int_{\mathbb{R}^{n}}J_{\varepsilon}(x)dx=1$. Since $J_{\varepsilon} * u_{\varepsilon}$ is a smooth function and the derivatives of each order are bounded, (\ref{3.1.5}) can be regarded as the advection equation (\ref{2.13}). We know that for any $u_0\in L^{1}(\mathbb{R}^{n})\cap H^{\gamma}(\mathbb{R}^{n})$, $u_0>0$, there is a unique global mild solution $u_{\varepsilon}$ in (\ref{3.1.5}). Further, this mild solution $u_{\varepsilon}$ is still a strong solution and for any $k\geq 0$, $u_{\varepsilon}\in C([0,T];C^{k}(\mathbb{R}^{n}))$.

Multiplying equation (\ref{3.1.5}) with $ru_{\varepsilon}^{r-1}\phi_{R}(x)$ and integrate over $\mathbb{R}^{n}$, we have
\begin{align}\label{3.1.6}
        &\int_{\mathbb{R}^{n}}ru_{\varepsilon}^{r-1}\partial_{t}^{\beta}u_{\varepsilon}\phi_{R}(x)dx\notag\\
        \leq{} & -\frac{4(r-1)}{r}\int_{\mathbb{R}^{n}}\lvert D^{\frac{\alpha}{2}}u_{\varepsilon}^{\frac{r}{2}}\rvert_{L^{2}}^{2}\phi_{R}(x)dx-r\int_{\mathbb{R}^{n}}u_{\varepsilon}^{r}D^{\alpha}\phi_{R}(x)dx\notag\\
        +{} & (r-1)\int_{\mathbb{R}^{n}}(J_{\varepsilon} * u_{\varepsilon})u_{\varepsilon}^{r}\phi_{R}(x)dx+ \int_{\mathbb{R}^{n}}n\nabla v_{\varepsilon}\cdot \nabla \phi_{R}(x)u_{\varepsilon}^{r}dx\notag\\
        +{}& ar \int_{\mathbb{R}^{n}}u_{\varepsilon}^{r}\phi_{R}(x)dx-br\int_{\mathbb{R}^{n}}u_{\varepsilon}^{r+1}\phi_{R}(x)dx\notag\\
        \leq{} & -\frac{4(r-1)}{r}\int_{\mathbb{R}^{n}}\lvert D^{\frac{\alpha}{2}}u_{\varepsilon}^{\frac{r}{2}}\rvert_{L^{2}}^{2}\phi_{R}(x)dx+(r-1)\int_{\mathbb{R}^{n}}(J_{\varepsilon} * u_{\varepsilon})u_{\varepsilon}^{r}\phi_{R}(x)dx\notag\\
        +{}& ar \int_{\mathbb{R}^{n}}u_{\varepsilon}^{r}\phi_{R}(x)dx-br\int_{\mathbb{R}^{n}}u_{\varepsilon}^{r+1}\phi_{R}(x)dx\notag\\
        +{}& \frac{rC}{R^{\alpha}}\Vert u \Vert_{L^{r}}+\frac{C}{R}\int_{\mathbb{R}^{n}}\lvert \nabla v_{\varepsilon} \rvert u_{\varepsilon}^{r}dx.
\end{align}
By using the Hardy-Littlewood-Sobolev inequality, we know
\begin{equation}
    \int_{\mathbb{R}^{n}}\lvert \nabla v_{\varepsilon}\rvert u_{\varepsilon}^{r}dx\leq \Vert u \Vert_{L^{\frac{nr}{n+1}}}^{r}.
\end{equation}

From $\Vert u_{\varepsilon}\Vert_{L^{r}}\leq C_{\varepsilon}$ and $\Vert u \Vert_{L^{\frac{nr}{n+1}}}\leq C_{\varepsilon}$, we know The last two terms in (\ref{3.1.6}) tend to 0 when $R\to 0$. Notice that $ \int_{\mathbb{R}^{n}}u_{\varepsilon}^{r}(J_{\varepsilon} * u_{\varepsilon})dx\leq \int_{\mathbb{R}^{n}}u_{\varepsilon}^{r+1}dx$. Thus (\ref{3.1.6}) can be written as
\begin{equation}
    \begin{split}
        \partial_{t}^{\beta}\int_{\mathbb{R}^{n}} u_{\varepsilon}^{r}dx+\frac{4(r-1)}{r}\int_{\mathbb{R}^{n}}\lvert D^{\frac{\alpha}{2}}u_{\varepsilon}^{\frac{r}{2}}\rvert_{L^{2}}^{2}dx &\leq (r-1-br)\int_{\mathbb{R}^{n}}u_{\varepsilon}^{r+1}dx\\
        &+ar\int_{\mathbb{R}^{n}}u_{\varepsilon}^{r}dx.
    \end{split}
\end{equation}
Under the assumption of $b\geq 1$, we have 
\begin{equation}
    \partial_{t}^{\beta}\int_{\mathbb{R}^{n}}u_{\varepsilon}^{r}dx+\frac{4 (r-1)}{r}\int_{\mathbb{R}^{n}}\lvert D^{\frac{\alpha}{2}}u_{\varepsilon}^{\frac{r}{2}}\rvert_{L^{2}}^{2}dx \leq ar\int_{\mathbb{R}^{n}}u_{\varepsilon}^{r}dx.
\end{equation}
For the initial data $u_{0}\in L^{1}(\mathbb{R}^{n})\cap L^{\infty}(\mathbb{R}^{n})$, the following estimates are obtained
\begin{equation}
    \Vert u_{\varepsilon}\Vert_{L^{\infty}(0,T;L^{1}(\mathbb{R}^{n})\cap L^{q}(\mathbb{R}^{n}))}\leq C,
\end{equation}
\begin{equation}
    \Vert u_{\varepsilon}\Vert_{L^{q}(0,T;L^{q}(\mathbb{R}^{n}))}\leq C,
\end{equation}
\begin{equation}
    \Vert D^{\frac{\alpha}{2}}u_{\varepsilon}^{\frac{r}{2}}\Vert_{L^{2}(0,T;L^{2}(\mathbb{R}^{n}))}\leq C,\ \ 1<r\leq q,
\end{equation}
\begin{equation}
    \Vert \partial_{t}^{\beta}u_{\varepsilon}\Vert_{L^{q_{1}}(0,T;W^{-\alpha,\frac{nq}{2n-q}}(\mathbb{R}^{n}))}\leq C, \ q_1>1.
\end{equation}
\textbf{Step 2.} (The application of the Aubin-Lions-Dubinski\u{\i} Lemma)\\
We first introduce $P_{+}\in B$, where $B$ is an Banach space, $B=L^{q}(\Omega)$, $\Omega$ is any bounded region, and define $P_{+}(\Omega):=\{u: [u]\leq C\}$ with $[u]=\Vert D^{\frac{\alpha}{2}}u^{\frac{q}{2}}\Vert_{L^{2}}^{\frac{q}{2}}+\Vert u \Vert_{L^{1}}+\Vert u \Vert_{L^{q}}$. The work in Huang \cite{huang2016well} has proved that $P_{+}(\Omega)\hookrightarrow\hookrightarrow L^{q}(\Omega)$. Recall that
\begin{equation}
    \Vert u_{\varepsilon} \Vert_{L^{q}(0,T;P_{+}(\Omega))}\leq C,
\end{equation}
\begin{equation}
     \Vert u_{\varepsilon} \Vert_{L^{q}(0,T;L^{q}(\Omega))}\leq C,
\end{equation}
\begin{equation}
    \Vert \partial_{t}^{\beta}u_{\varepsilon} \Vert_{L^{q_1}(0,T;W^{-\alpha, \frac{nq}{2n-q}}(\Omega))}\leq C, \ q_{1}>1.
\end{equation}
By Sobolv's embedding theorem, there is $P_{+}(\Omega)\hookrightarrow\hookrightarrow L^{q}(\Omega)\hookrightarrow\hookrightarrow W^{-\alpha, \frac{nq}{2n-q}}(\Omega)$ and by the Aubin-Lions-Dubinski\u{\i} Lemma, there exists a subsequence $u_{\varepsilon}$ such that 
$$u_{\varepsilon} \to u\ \ {\rm in}\  L^{q}(0,T;L^{q}(\Omega)),$$
as $\varepsilon \to 0$. Let $\{B_{k}\}_{k=1}^{\infty} \in \mathbb{R}^{n}$ be a sequence of balls centered 0 with radius $R_{k},$ $R_{k}\to \infty$, there exists a subsequence $u_{\varepsilon}$ satisfies the following strong convergence
$$u_{\varepsilon} \to u \ \ {\rm{in}} \ L^{q}(0,T;L^{(B_{k})}), \ \ {\rm{as}}\  \varepsilon\to 0,\ \forall k.$$
\textbf{Step 3.} (The existence of a global weak solution)\\
Next, we will prove that $u$ is a weak solution to (\ref{1.1}). For any $\phi \in C_{c}^{\infty}(\mathbb{R}^{n})$ and any $0<t<\infty$
\begin{equation}\label{3.2.8}
\begin{split}
    \int_{0}^{T}\int_{\mathbb{R}^{n}}(u_{\varepsilon}-u_{0})\tilde{\partial}_{T}^{\beta} \phi dxdt
    &= - \int_{0}^{T}\int_{\mathbb{R}^{n}}[u_{\varepsilon}D^{\alpha}\phi(x)]dxdt\\
    &+\int_{0}^{T}\int_{\mathbb{R}^{n}}u_{\varepsilon}\nabla v_{\varepsilon} \cdot \nabla \phi(x)dxdt\\
    & + a\int_{0}^{T}\int_{\mathbb{R}^{n}}u_{\varepsilon}\phi(x)dxdt-b\int_{0}^{T}\int_{\mathbb{R}^{n}}u_{\varepsilon}^{2}\phi(x)dxdt.
\end{split}
\end{equation}
(i) $$\int_{\mathbb{R}^{n}}[u_{\varepsilon}D^{\alpha}\phi(x)]dx\to \int_{\mathbb{R}^{n}}[uD^{\alpha}\phi(x)]dx,$$
as $\varepsilon \to 0.$\\
(ii) $$\int_{\mathbb{R}^{n}}u_{\varepsilon}\nabla v_{\varepsilon} \cdot \nabla \phi(x)dx\to \int_{\mathbb{R}^{n}}u\nabla v \cdot \nabla \phi(x)dx,$$
as $\varepsilon \to 0.$\\
The step (i) and (ii) follow the proof of Theorem 2.3 in \cite{huang2016well}. Since $u_{\varepsilon}\to u$ in $L^{q}(0,T;L^{q}(\Omega))$\\
(iii)
\begin{equation}
\begin{split}
\left\lvert a\int_{\mathbb{R}^{n}}u_{\varepsilon}\phi(x)dx-a\int_{\mathbb{R}^{n}}u\phi(x)dx\right\rvert&\leq a\int_{\mathbb{R}^{n}}\lvert u_{\varepsilon}-u \rvert\phi(x)dx\\
&\leq C\int_{\Omega}\lvert u_{\varepsilon}-u\rvert dx\\
&\leq C(\Omega)\Vert u_{\varepsilon}-u\Vert_{L^{q}(0,T;L^{q}(\Omega))}\\
&\to 0,
\end{split}
\end{equation}
as $\varepsilon\to 0.$\\
(iv) 
\begin{equation}
   \begin{split}
       \left\lvert b\int_{\mathbb{R}^{n}}u_{\varepsilon}^{2}\phi(x)dx-b\int_{\mathbb{R}^{n}}u^{2}\phi(x)dx\right\rvert&\leq b\int_{\mathbb{R}^{n}}\lvert u_{\varepsilon}^{2}-u^{2} \rvert\phi(x)dx\\
       &\leq C\int_{\Omega}\lvert u_{\varepsilon}-u\rvert (u_{\varepsilon}+u)dx\\
       &\leq C(\Omega)\Vert u_{\varepsilon}-u\Vert_{L^{q}(0,T;L^{q}(\Omega))}\\
       &\to 0,
   \end{split}
\end{equation}
as $\varepsilon\to 0.$\\
Then combine (i)-(iv) and Let $\varepsilon\to 0$, for any $0<t<\infty$, we have 
\begin{equation}
\begin{split}
    \int_{0}^{T}\int_{\mathbb{R}^{n}}(u-u_{0})\tilde{\partial}_{t}^{\beta} \phi dxdt
    &= - \int_{0}^{T}\int_{\mathbb{R}^{n}}[uD^{\alpha}\phi(x)]dxdt\\
    &+\int_{0}^{T}\int_{\mathbb{R}^{n}}u\nabla v \cdot \nabla \phi(x)dxdt\\
    & + a\int_{0}^{T}\int_{\mathbb{R}^{n}}u\phi(x)dxdt-b\int_{0}^{T}\int_{\mathbb{R}^{n}}u^{2}\phi(x)dxdt.
    \end{split}
\end{equation}
Which means that $(u,v)$ is a global weak solution.\\
\textbf{Step 4.} (The uniformly in time $L^{\infty}$ estimate of weak solution)\\
We denote 
$$q_{k}:=2^{k}+\max\left\{\frac{n}{\alpha}, \frac{1}{1-b}\right\}.$$
Multiplying the first equation in (\ref{1.1}) by $q_{k}u^{q_{k}-1}$, we have 
\begin{equation}
    \begin{split}
        \int_{\mathbb{R}^{n}}q_{k}u^{q_{k}-1} \partial_{t}^{\beta} udx&= -\frac{4(q_{k}-1)}{q_{k}}\int_{\mathbb{R}^{n}}\lvert D^{\frac{\alpha}{2}}u^{\frac{q_{k}}{2}}\rvert^{2}dx\\
        &+(q_{k}-1-bq_{k})\int_{\mathbb{R}^{n}}u^{q_{k}+1}dx+aq_{k}\int_{\mathbb{R}^{n}}u^{q_{k}}dx.
    \end{split}
\end{equation}
$q_{k}-1-bq_{k}<0$ from $b\geq 1$, by Lemma \ref{2.10}, one has
\begin{equation}\label{3.1.7}
    \begin{split}
        \partial_{t}^{\beta}\int_{\mathbb{R}^{n}}u^{q_{k}}dx &\leq -\frac{4(q_{k}-1)}{q_{k}}\int_{\mathbb{R}^{n}}\lvert D^{\frac{\alpha}{2}}u^{\frac{q_{k}}{2}}\rvert^{2}dx\\
        &+aq_{k}\int_{\mathbb{R}^{n}}u^{q_{k}}dx.
    \end{split}
\end{equation}
(\ref{3.1.7}) can be rewritten as
\begin{equation}\label{3.1.8}
\begin{split}
_{0}^{c}\textrm{D}_{t}^{\beta}\Vert u \Vert_{L^{q_{k}}(\mathbb{R}^{n})}^{q_{k}} &\leq -\frac{4(q_{k}-1)}{q_{k}}\Vert  D^{\frac{\alpha}{2}}u^{\frac{q_{k}}{2}} \Vert_{L^{2}(\mathbb{R}^{n})}^{2}+aq_{k}\Vert u \Vert_{L^{q_{k}}(\mathbb{R}^{n})}^{q_{k}}\\
&\leq -2C_{1}\Vert  D^{\frac{\alpha}{2}}u^{\frac{q_{k}}{2}} \Vert_{L^{2}(\mathbb{R}^{n})}^{2}+aq_{k}\Vert u \Vert_{L^{q_{k}}(\mathbb{R}^{n})}^{q_{k}},
\end{split}
\end{equation}
where $0<C_{1}\leq \frac{2(q_{k}-1)}{q_{k}}$ is a fixed constant. Using the Interpolation inequality and Sobolev inequality to the last term $\Vert u \Vert_{L^{q_{k}}(\mathbb{R}^{n})}^{q_{k}}$.
\begin{equation}\label{3.1.10}
   \begin{split}
       \Vert u \Vert_{L^{q_{k}}(\mathbb{R}^{n})}^{q_{k}}&\leq \Vert u \Vert_{L^{\frac{nq_{k}}{n-\alpha}}(\mathbb{R}^{n})}^{q_{k}\tau_{1}}\Vert u \Vert_{L^{q_{k-1}}(\mathbb{R}^{n})}^{q_{k}(1-\tau_1)}\\
       & =\Vert u^{\frac{q_{k}}{2}} \Vert_{L^{\frac{2n}{n-\alpha}}(\mathbb{R}^{n})}^{2\tau_{1}}\Vert u \Vert_{L^{q_{k-1}}(\mathbb{R}^{n})}^{q_{k}(1-\tau_1)}\\
       & \leq S_{\alpha,n}^{2\tau_{1}}\Vert D^{\frac{\alpha}{2}}u^{\frac{q_{k}}{2}}\Vert_{L^{2}(\mathbb{R}^{n})}^{2\tau_{1}}\Vert u\Vert_{L^{q_{k-1}}(\mathbb{R}^{n})}^{q_{k}(1-\tau_1)},
   \end{split} 
\end{equation}
where $$\tau_{1}=\frac{\frac{1}{q_{k-1}}-\frac{1}{q_{k}}}{\frac{1}{q_{k-1}}-\frac{n-\alpha}{nq_{k}}}.$$
The Young's inequality tells that
\begin{equation}\label{3.1.9}
    \begin{split}
        aq_{k}\Vert u\Vert_{L^{q_{k}}(\mathbb{R}^{n})}^{q_{k}}&\leq\frac{1}{c_{1}}\theta_{1}^{c_{1}}\Vert D^{\frac{\alpha}{2}}u^{\frac{q_{k}}{2}}\Vert_{L^{2}(\mathbb{R}^{n})}^{2c_{1}\tau_{1}}+\frac{1}{c_{2}}\theta_{1}^{-c_2}S_{\alpha,n}^{2\tau_{1}\theta_{1}c_{2}}(aq_{k})^{c_2}\Vert u\Vert_{L^{q_{k-1}}(\mathbb{R}^{n})}^{c_{2}q_{k}(1-\tau_1)}\\
        &\leq C_{1}\Vert D^{\frac{\alpha}{2}}u^{\frac{q_{k}}{2}}\Vert_{L^{2}(\mathbb{R}^{n})}^{2}+C_{2}(q_{k})(aq_{k})^{c_2}\Vert u\Vert_{L^{q_{k-1}}(\mathbb{R}^{n})}^{c_{2}q_{k}(1-\tau_1)},
    \end{split}
\end{equation}
where $c_{1}=\frac{1}{\tau_1}$ and $$c_2=\frac{1}{1-\tau_1}=\frac{2\alpha q_{k-1}}{2nq_{k}-2(n-\alpha)q_{k-1}}\leq 1+n$$
$$\theta_{1}=(C_{1}c_{1})^{\frac{1}{c_{1}}},\ \ C_{2}(q_{k})=\frac{1}{c_{2}}(C_{2}c_{2})^{-\frac{c_2}{c_{1}}}S_{\alpha,n}^{2\tau_{1}\theta_{1}c_2}.$$
We see that $C_{2}(q_{k})$ is uniformly bounded since $c_{1}\to\frac{n+1}{n}$ and $c_{2}\to n+1$ as $k\to \infty.$ Substituting (
\ref{3.1.9}) into (\ref{3.1.8}) yields to
\begin{equation}
    _{0}^{c}\textrm{D}_{t}^{\beta}\Vert u\Vert_{L^{q_{k}}(\mathbb{R}^{n})}^{q_{k}} \leq -C_{1}\Vert D^{\frac{\alpha}{2}}u^{\frac{q_{k}}{2}}\Vert_{L^{2}(\mathbb{R}^{n})}^{2}+C_{2}(q_{k})(aq_{k})^{c_2}\left(\Vert u \Vert_{L^{q_{k-1}}(\mathbb{R}^{n})}^{q_{k-1}}\right)^{\eta_1},
\end{equation}
where 
$$\eta_{1}=\frac{q_{k}(1-\tau_{1})c_{2}}{q_{k-1}}=\frac{q_{k}}{q_{k-1}}<2.$$
It can be seen from (\ref{3.1.9}) that
\begin{equation}
    \Vert u\Vert_{L^{q_{k}}(\mathbb{R}^{n})}^{q_{k}}\leq C_{1}\Vert D^{\frac{\alpha}{2}}u^{\frac{q_{k}}{2}}\Vert_{L^{2}(\mathbb{R}^{n})}^{2}+C_{2}(q_{k})\left(\Vert u \Vert_{L_{q^{k-1}}(\mathbb{R}^{n})}^{q_{k-1}}\right)^{\eta_1}.
\end{equation}
We have
\begin{equation}
\begin{split}
    _{0}^{c}\textrm{D}_{t}^{\beta} \Vert u\Vert_{L^{q_{k}}(\mathbb{R}^{n})}^{q_{k}}&\leq - \Vert u\Vert_{L^{q_{k}}(\mathbb{R}^{n})}^{q_{k}}+C_{2}(q_{k})\left((aq_{k})^{c_2}+1\right)\left(\Vert u \Vert_{L^{q_{k-1}}(\mathbb{R}^{n})}^{q_{k-1}}\right)^{\eta_1}\\
    &\leq - \Vert u\Vert_{L^{q_{k}}(\mathbb{R}^{n})}^{q_{k}}+C_{2}(q_{k})\left(a^{c_2}+1\right)q_{k}^{c_2}\left(\Vert u \Vert_{L^{q_{k-1}}(\mathbb{R}^{n})}^{q_{k-1}}\right)^{\eta_1}.
\end{split}
\end{equation}
Let $C(K)=c_{2}(q_{k})(a^{c_2}+1)$, $C(K)\to C_{0}$, $k\to \infty$, where $C_{0}$ is a constant independent of $q_{k}$. Therefore, we can see that $C(K)$ has a uniform upper bound $C$ for $k>1$, that is
\begin{equation}
    _{0}^{c}\textrm{D}_{t}^{\beta}\Vert u\Vert_{L^{q_{k}}(\mathbb{R}^{n})}^{q_{k}}\leq -\Vert u\Vert_{L^{q_{k}}(\mathbb{R}^{n})}^{q_{k}}+Cq_{k}^{c_2}\left(\Vert u \Vert_{L^{q_{k-1}}(\mathbb{R}^{n})}^{q_{k-1}}\right)^{\eta_1}.
\end{equation}
\textbf{Step 5.} (The uniform $L^{\infty}$ estimate)\\
Let$$y_{k}(t) :=\Vert u\Vert_{L^{q_{k}}(\mathbb{R}^{n})}^{q_{k}}.$$
We have
\begin{equation}
     _{0}^{c}\textrm{D}_{t}^{\beta}y_{k}(t)+y_{k}(t)\leq Cq_{k}^{c_2}y_{k-1}^{\eta_1}(t).
\end{equation}
By Lemma \ref{2.12}, we obtain
\begin{equation}
    \begin{split}
        y_{k}(t) &\leq y_{k}(0)+\frac{1}{\Gamma(\beta)}\int_{0}^{t}(t-s)^{\beta-1}Cq_{k}^{c_2}y_{k-1}^{\eta_1}(s)ds\\
        &\leq y_{k}(0)+Cq_{k}^{c_2}\max\{\underset{t\geq 0}{\sup} y_{k-1}^{2}(t),1\}\\
        &\leq Cq_{k}^{c_2}\max\{1,y_{k}(0),\underset{t\geq 0}{\sup} y_{k-1}^{2}(t)\}.\\
        &\leq Cq_{k}^{1+n}\max\{1,y_{k}(0),\underset{t\geq 0}{\sup} y_{k-1}^{2}(t)\}.
        \end{split}
    \end{equation}
 Let $\mu :=\max\left\{\frac{\alpha}{n}, \frac{1}{1-b}\right\}+1$, then $q_{k}\leq 2^{k}\mu$, one has
 \begin{equation}
     y_{k}(t)\leq C\left(2^{k}\mu\right)^{n+1}\max\{1 ,y_{k}(0),\underset{t\geq 0}{\sup} y_{k-1}^{2}(t)\}
 \end{equation}
 Take $a_{k}:=C\left(2^{k}\mu\right)^{n+1}$, denote $K_{0} := \max\{1, \Vert u_{0}\Vert_{L^{1}(\mathbb{R}^{n})},  \Vert u_{0}\Vert_{L^{1}(\mathbb{R}^{n})}\}$, $y_{k}(0)=\Vert u_{0}\Vert_{L^{q_{k}}}^{q_{k}}\leq \Vert u_{0} \Vert_{L^{1}} \Vert u_{0}\Vert_{L^{\infty}}^{q_{k}-1}$, we have
 $$\max\{y_{k}(0),1\}\leq A^{q_{k}}.$$
 So we obtain
 \begin{equation}
     \begin{split}
         y_{k}(t)&\leq a_{k}\max\{\underset{t\geq 0}{\sup} y_{k-1}^{2}(t), A^{q_{k}}\}\\
         &\leq (C\mu^{n+1})^{2^{k}-1}2^{(n+1)(2^{k}-k-2)}\max\{\underset{t\geq 0}{\sup} y_{0}^{2^k}(t),\sum_{i=0}^{k-1}A^{q_{k-i}2^{i}}\}
     \end{split}
 \end{equation}
 Since $A^{q_{k-i}2^{i}}\leq\widetilde{A}^{q_{k}}$, where $\widetilde{A}>1$ is constant independent of $k$. We have
 \begin{equation}
     \Vert u\Vert_{L^{q_{k}}}^{q_{k}}\leq (c\mu^{n+1})^{2^{k}-1}2^{(n+1)(2^{k+1}-k-2)}\max\{\underset{t\geq 0}{\sup} y_{0}^{2^k}(t),K\widetilde{A}^{q_{k}}\}.
 \end{equation}
 Taking both sides to the power of $\frac{1}{q_{k}}$ at the same time, we have
 \begin{equation}\label{3.10}
     \begin{split}
             \Vert u \Vert_{L^{q_{k}}}&\leq (C\mu^{n+1})^{\frac{2^{k}-1}{q_{k}}}2^{\frac{(n+1)(2^{k}-k-2)}{q_{k}}}\max\{\underset{t\geq 0}{\sup}y_{0}^{\frac{2^{k+1}}{q_{k}}}, K^{\frac{1}{q_{k}}}\widetilde{A}\}\\
             &\leq C\mu^{n+1}2^{2(n+1)}\max\{\underset{t\geq 0}{\sup}y_{0}(t), K^{\frac{1}{q_{k}}}\widetilde{A}\} 
     \end{split}
 \end{equation}
 On the other hand, 
 \begin{equation}\label{3.11}
 y_{0}(t)=\int_{\mathbb{R}^{n}}u^{q_0}dx=\int_{\mathbb{R}^{n}}u^{\beta}dx\leq C.
 \end{equation}
 From (\ref{3.10}) and (\ref{3.11}), we know that there is a constant $C$ independent of $q_{k}$ such that $$\Vert u \Vert_{L^{q_{k}}}\leq C.$$
 Therefore, when $k\to \infty$, there is $u\in L^{\infty}(0,T; L^{\infty}(\mathbb{R}^{n}))$. That is,
 $$\Vert u \Vert_{L^{\infty}(0,T;L^{\infty}(\mathbb{R}^{n}))}\leq C.$$\qed

\subsection{The case of \ \texorpdfstring {$1-\frac{\alpha}{n}<b<1$}{lg}}
\begin{lemma}\label{3.2}
Under the assumption of Theorem \ref{1.1.1}(ii), the solution $(u,v)$ of system (\ref{1.1}) satisfies the following estimation
\begin{flalign}\label{3.2.1}
        && \Vert u \Vert_{L^{\infty}(0,T;L^{q}(\mathbb{R}^{n}))} &\leq C, \ \ \Vert D^{\frac{\alpha}{2}}u^{\frac{q}{2}} \Vert_{L^{2}(0,T;L^{2}(\mathbb{R}^{n}))}\leq C.
    \end{flalign}
    \begin{equation}\label{3.2.2}
        \Vert \nabla v \Vert_{L^{\infty}([0,\infty];L^{s_1}(\mathbb{R}^{n}))}\leq C, \ \ s_1 > \frac{n}{n-1}.
    \end{equation}
\end{lemma}
\begin{pro}
We divide into two cases:
    \begin{enumerate}[(1)]
        \item $1 < q \leq \frac{1}{1-b}$
        
        In this case, we know $q-1-bq<0$, Then the proof process is similar to Lemma \ref{3.1}, that is, (\ref{3.2.1}) holds.
        \item $q=\frac{n}{\alpha}>\frac{1}{1-b}$
         
         From $1-\frac{\alpha}{n}<b<1$ we know $1<\frac{n}{\alpha}<\frac{1}{1-b}$, therefore, $$\Vert u \Vert_{L^{\infty}(0,T;L^{\frac{n}{\alpha}}(\mathbb{R}^{n}))} \leq C.$$ Taking a sufficiently small $\varepsilon>0$,such that $1<\frac{n}{\alpha}+\varepsilon<\frac{1}{1-b},$ we obtain
         \begin{equation}\label{3.2.3}
            \Vert u \Vert_{L^{\infty}(0,T;L^{\frac{n}{\alpha}+\epsilon}(\mathbb{R}^{n}))} \leq C.
        \end{equation}
    Next using Interpolation inequality and Sobolev inequality to $\Vert u \Vert_{L^{q+1}(\mathbb{R}^{n})}$
    \begin{equation}\label{3.2.4}
        \begin{split}
            \Vert u \Vert_{L^{q+1}}^{q+1} & \leq \Vert u \Vert_{L^{\frac{nq}{n-\alpha}}}^{(q+1)\tau_2}\Vert u \Vert_{L^{\frac{n}{\alpha}+\epsilon}}^{(q+1)(1-\tau_2)}\\
            &= \Vert u^{\frac{q}{2}}\Vert_{L^{\frac{2n}{n-\alpha}}}^\frac{2(q+1)\tau_2}{q} \Vert u \Vert_{L^{\frac{n}{\alpha}+\epsilon}}^{(q+1)(1-\tau_2)}\\
            &\leq S_{\alpha, n}^{\frac{2(q+1)\tau_2}{q}}\Vert D^{\frac{\alpha}{2}} u^{\frac{q}{2}}\Vert_{L^{2}}^\frac{2\tau_2(q+1)}{q} \Vert u \Vert_{L^{\frac{n}{\alpha}+\epsilon}}^{(q+1)(1-\tau_2)},
        \end{split}
\end{equation}
where $\tau_2$ satisfies $\frac{1}{q+1}=\frac{(n-\alpha)\tau_2}{nq}+\frac{1-\tau_2}{\frac{n}{\alpha}+\epsilon}$, it is easy to see that $$\frac{2(q+1)\tau_2}{q}=\frac{\frac{4(q+1)}{n+\alpha\epsilon}-2}{\frac{2q}{n+\alpha\epsilon}-\frac{n-\alpha}{n}}<2.$$
Then using Young's inequality, we have
\begin{equation}\label{3.2.5}
       \begin{split}
           (q-1-bq)\Vert u \Vert_{L^{q+1}}^{q+1} 
           &\leq \frac{1}{c_3}\theta_{2}^{c_3}\Vert D^{\frac{\alpha}{2}}u^{\frac{q}{2}}\Vert_{L^{2}}^{2(q+1)\tau_{2}c_{3}}\\
           &+\frac{1}{c_4}\theta_{2}^{-c_4}S_{\alpha, n}^{\frac{2(q+1)\tau_2c_4}{q}}(q-1-bq)^{c_4}\Vert u \Vert_{L^{\frac{n}{\alpha}+\epsilon}(\mathbb{R}^{n})}^{qc_4(1-\tau_2)}\\
           & \leq C_3\Vert D^{\frac{\alpha}{2}}u^{\frac{q}{2}}\Vert_{L^{2}}^{2(q+1)\tau_2c_3}+C_{4}(q-1-bq)^{c_4}\Vert u \Vert_{L^{\frac{n}{\alpha}+\epsilon}}^{(q+1)(1-\tau_2)c_4},
       \end{split}
   \end{equation}
   where $c_3, c_4>1$ satisfies $\frac{1}{c_3}+\frac{1}{c_4}=1$. Take $c_3=\frac{q}{\tau_2(q+1)}$, $C_3=\frac{2(p-1)}{q}$. We all know
   \begin{equation}\label{3.2.6}
       \begin{split}
           \int_{\mathbb{R}^{n}}qu^{q-1}\partial_{t}^{\beta}udx+\frac{4(q-1)}{q}\int_{\mathbb{R}^{n}}\left\lvert D^{\frac{\alpha}{2}}u^{\frac{q}{2}}\right\rvert^{2}dx&=(q-1-bq)\int_{\mathbb{R}^{n}}u^{q+1}dx\\
           &+aq\int_{\mathbb{R}^{n}}u^{q}dx.
       \end{split}
   \end{equation}
 By Lemma \ref{2.10} and substituting (\ref{3.2.5}) into (\ref{3.2.6}) yields to
 \begin{equation}\label{3.2.6.1}
     \partial_{t}^{\beta}\int_{\mathbb{R}^{n}}u^{q}dx+\frac{2(q-1)}{q}\int_{\mathbb{R}^{n}}\left\lvert D^{\frac{\alpha}{2}}u^{\frac{q}{2}}\right\rvert^{2}dx\leq C+aq\int_{\mathbb{R}^{n}}u^{q}dx.
 \end{equation}
 The above formula contains
 $$\partial_{t}^{\beta}\int_{\mathbb{R}^{n}}u^{q}dx\leq C+aq\int_{\mathbb{R}^{n}}u^{q}dx.$$
 Also by Lemma \ref{2.11}, we obtain $\Vert u \Vert_{L^{\infty}(0,T;L^{q}(\mathbb{R}^{n}))}\leq C.$ It is easy to see that $\left\Vert D^{\frac{\alpha}{2}}u^{\frac{q}{2}}\right\Vert_{L^{2}(0,T;L^{2}(\mathbb{R}^{n}))}\leq C.$
 
 The following proof process is similar to Lemma \ref{3.1}, so we can get that (\ref{3.2.2}) holds.\qed
 \end{enumerate}
\end{pro}
Regarding the proof process for the estimation of the two necessary conditions for applying the Aubin-Lions-Dubinski\u{\i} Lemma, the spatial derivative and the time derivative, similar to Lemma \ref{3.2.1.1.2} can be obtained

\noindent $\mathbf{Proof\ of\ Theorem\ref{1.1.1}(ii):}$\\
\textbf{Step 1.} (Regularization and a priori estimates)\\ We also consider the regularization problem (\ref{3.1.5}), combine $1-\frac{\alpha}{n}<b<1$, we obtain
\begin{equation}
    \partial_{t}^{\beta}\Vert u_{\varepsilon}\Vert_{L^{q}}^{q}+\frac{4(q-1)}{q}\Vert D^{\frac{\alpha}{2}}u_{\varepsilon}^{\frac{q}{2}}\Vert_{L^{2}}^{2}\leq (q-1-bq)\Vert u_{\varepsilon}\Vert_{L^{q+1}}^{q+1}+aq\Vert u_{\varepsilon}\Vert_{L^{q}}^{q}.
\end{equation}
The following proof process is the same as the proof method of the second case of Lemma \ref{3.2}, it can be obtained that the solution $(u_{\varepsilon},v_{\varepsilon})$ of the regularization problem satisfies
\begin{equation}
    \Vert u_{\varepsilon}\Vert_{L^{\infty}(0,T;L^{1}(\mathbb{R}^{n})\cap L^{q}(\mathbb{R}^{n}))}\leq C,
\end{equation}
\begin{equation}
    \Vert D^{\frac{\alpha}{2}}u_{\varepsilon}^{\frac{q}{2}}\Vert_{L^{2}(0,T;L^{2}(\mathbb{R}^{n}))}\leq C,
\end{equation}
\begin{equation}
    \Vert \partial_{t}^{\beta}u_{\varepsilon}\Vert_{L^{q_1}(0,T;W^{-\alpha,\frac{nq}{2n-q}}(\mathbb{R}^{n}))}\leq C,\ q_1>1.
\end{equation}
\textbf{Step 2.} (The application of the Aubin-Lions-Dubinski\u{\i} Lemma)\\
Here is the same as the proof process of Theorem \ref{1.1.1}(i), by 
\begin{equation}
    \Vert u_{\varepsilon} \Vert_{L^{q}(0,T;P_{+}(\Omega))}\leq C,
\end{equation}
\begin{equation}
     \Vert u_{\varepsilon} \Vert_{L^{q}(0,T;L^{q}(\Omega))}\leq C,
\end{equation}
\begin{equation}
    \Vert \partial_{t}^{\beta}u_{\varepsilon}\Vert_{L^{q_1}(0,T;W^{-\alpha,\frac{nq}{2n-q}}(\mathbb{R}^{n}))}\leq C,\ q_1>1.
\end{equation}
and $P_{+}(\Omega)\hookrightarrow\hookrightarrow L^{q}(\Omega)\hookrightarrow\hookrightarrow W^{-\alpha, \frac{nq)}{2n-q}}(\Omega)$, also by Aubin-Lions-Dubinski\v i Lemma, we know
\begin{equation}
    u_{\varepsilon} \to u \ \ {\rm{in}} \ L^{q}(0,T;L^{q}(\Omega)), \ \ {\rm{as}}\  \varepsilon\to 0.
\end{equation}
Similarly, there exists a subsequence $u_{\varepsilon}$ satisfies the following strong convergence
$$u_{\varepsilon} \to u \ \ {\rm{in}} \ L^{q}(0,T;L^{q}(B_{k})), \ \ {\rm{as}}\  \varepsilon\to 0.$$
\textbf{Step 3.} (The existence of a global weak solution)\\
The weak formulation for $u_{\varepsilon}$ is that for any test formulation $\varphi(x) \in C_{c}^{\infty}(\mathbb{R}^{n})$ and $0<T<\infty$, we have (\ref{3.2.8}) holds. The remaining proof steps are the same as the proof process of Theorem \ref{1.1.1} (i).\\
\textbf{Step 4.} (The uniformly in time $L^{\infty}$ estimate of weak solution)\\
Take $q$ in  (\ref{3.1.3}) as $q_k$, also use interpolation inequality and Sobolev inequality to estimate $\Vert u \Vert_{L^{q_{k}+1}}^{q_{k}+1}$
\begin{equation}
    \begin{split}
        \Vert u \Vert_{L^{q_{k}+1}}^{q_{k}+1}&\leq \left\Vert u^{\frac{q_{k}}{2}} \right\Vert_{L^{\frac{2n}{n-\alpha}}}^{\frac{2(q_{k}+1)\tau_3}{q_{k}}}\Vert u \Vert_{L^{q_{k-1}}}^{(1-\tau_3)(q_{k}+1)}\\
        &\leq S_{\alpha,n}^{\frac{2(q_{k}+1)\tau_3}{q_{k}}}\left\Vert D^{\frac{\alpha}{2}}u^{\frac{q_{k}}{2}}\right\Vert_{L^{2}}^{\frac{2(q_{k}+1)\tau_3}{q_{k}}}\Vert u \Vert_{L^{q_{k-1}}}^{(1-\tau_3)(q_{k}+1)}.
    \end{split}
\end{equation}
where $$\tau_3=\frac{\frac{1}{q_{k-1}}-\frac{1}{q_{k}+1}}{\frac{1}{q_{k-1}}-\frac{n-\alpha}{n}}\sim O(1),\ \ 1-\tau_3\sim O(1), q_{k}\to \infty.$$
The Young's inequality tells that
\begin{equation}\label{3.2.9}
\begin{split}
    (q_{k}-1-bq_{k}) \Vert u \Vert_{L^{q_{k}+1}}^{q_{k}+1}&\leq \frac{1}{c_5}\theta_{3}^{c_5}\left\Vert D^{\frac{\alpha}{2}}u^{\frac{q_{k}}{2}}\right\Vert_{L^{2}}^{\frac{2c_5(q_{k}+1)\tau_3}{q_{k}}}\\
    &+\frac{1}{c_6}\theta_{3}^{-c_6}S_{\alpha, n}^{\frac{2c_6\tau_3(q_{k}+1)}{q_{k}}}\Vert u \Vert_{L^{q_{k-1}}}^{c_6(1-\tau_3)(q_{k}+1)},
\end{split}
\end{equation}
where $\frac{2c_5(q_{k}+1)\tau_3}{q_{k}}=2$ is taken, i.e., $c_5=\frac{q_{k}}{\tau_3(q_{k}+1)}$ and
$$c_6=\frac{q_k}{q_{k}-\tau_3(q_k+1)}=\frac{nq_{k}-(n-\alpha)q_{k-1}}{2q_{k-1}-n}\leq 1+n.$$
Substituting (\ref{3.2.9}) into (\ref{3.1.3}) yields to
\begin{equation}
    {_{0}^{c}\textrm{D}_{t}^{\beta}}\Vert u \Vert_{L^{q_{k}}}^{q_{k}}\leq -C_5\left\Vert D^{\frac{\alpha}{2}}u^{\frac{q_{k}}{2}}\right\Vert_{L^{2}}^{2}+C_6q_{k}^{c_6}\left(\Vert u \Vert_{L^{q_{k-1}}}^{q_{k-1}}\right)^{\eta_1}+aq_{k}\Vert u \Vert_{L^{q_{k}}}^{q_{k}},
\end{equation}
where $$\eta_1=\frac{c_6(1-\tau_3)(q_{k}+1)}{q_{k-1}}=\frac{q_{k}(1-\tau_3)}{q_{k}(1-\tau_3)-\tau_3}\frac{q_{k}+1}{q_{k-1}}<2.$$
From (\ref{3.1.10}) and (\ref{3.1.9}), we can further obtain
\begin{equation}
{_{0}^{c}\textrm{D}_{t}^{\beta}}\Vert u \Vert_{L^{q_{k}}}^{q_{k}}\leq -\Vert u \Vert_{L^{q_{k}}}^{q_{k}}+C_2q_{k}^{c_2}\left(\Vert u \Vert_{L^{q_{k-1}}}^{q_{k-1}}\right)^{\eta_1}+C_6q_{k}^{c_6}\left(\Vert u \Vert_{L^{q_{k-1}}}^{q_{k-1}}\right)^{\eta_2}.
\end{equation}
Since $C_2, C_6$ are uniformly bounded when $k\to \infty$. Take $C=\max\{C_2,C_6\}$, $l=\max\{c_2, c_6\}$, there is
\begin{equation}
    {_{0}^{c}\textrm{D}_{t}^{\beta}}\Vert u \Vert_{L^{q_{k}}}^{q_{k}}\leq -\Vert u \Vert_{L^{q_{k}}}^{q_{k}}+Cq_{k}^{l}\left[\left(\Vert u \Vert_{L^{q_{k-1}}}^{q_{k-1}}\right)^{\eta_1}+\left(\Vert u \Vert_{L^{q_{k-1}}}^{q_{k-1}}\right)^{\eta_2}\right].
\end{equation}
\textbf{Step 5.} (The uniform $L^{\infty}$ estimate)\\
Denote $y_{k}(t) :=\Vert u \Vert_{L^{q_{k}}}^{q_{k}}$, we have
\begin{equation}
     {_{0}^{c}\textrm{D}_{t}^{\beta}}y_{k}(t)+y_{k}(t)\leq Cq_{k}^{l}\left(y_{k-1}^{\eta_1}(t)+y_{k-1}^{\eta_2}(t)\right).
\end{equation}
By Lemma \ref{2.12}, we have
\begin{equation}
    \begin{split}
        y_{k}(t) &\leq y_{k}(0)+\frac{Cq_{k}^{l}}{\Gamma(\beta)}\int_{0}^{t}(t-s)^{\beta-1}\left(y_{k-1}^{\eta_1}(s)+y_{k-1}^{\eta_2}(s)\right)ds\\
        &\leq y_{k}(0)+2Cq_{k}^{l}\max\{1,\underset{t\geq 0}{\sup}y_{k-1}^{2}(t)\}\\
        &\leq 2Cq_{k}^{l}\max\{1, y_{k}(0), \underset{t\geq 0}{\sup}y_{k-1}^{2}(t)\}.
\end{split}
\end{equation}
Let $\mu :=\max\{\frac{n}{\alpha},\frac{1}{1-b}\}+1$, then $q_{k}\leq 2^{k}\mu$,
\begin{equation}
    y_{k}(t)\leq 2C(2^{k}\mu)^{n+1}\max\{1, y_{k}(0), \underset{t\geq 0}{\sup}y_{k-1}^{2}(t)\}.
\end{equation}
Let $a:=2C(2^{k}\mu)^{n+1}$, define $K_0 :=\max\{1, \Vert u_{0}\Vert_{L^{1}(\mathbb{R}^{n})}, \Vert u_{0} \Vert_{L^{\infty}(\mathbb{R}^{n})}\}$, $y_{k}(0)=\Vert u_{0}\Vert_{L^{q_{k}}}^{q_{k}}\leq \Vert u_0\Vert_{L^{1}}\Vert u_0 \Vert_{L^{\infty}}^{q_{k}-1}$, so we have $\max\{y_{k}(0), 1\}\leq K^{q_{k}}$, where constant $K>1$ is independent of $k$ but depends on $\Vert u_{0}\Vert_{L^{1}(\mathbb{R}^{n})}, \Vert u_{0} \Vert_{L^{\infty}(\mathbb{R}^{n})}$, therefore there is
\begin{equation}
    \begin{split}
        y_{k}(t)&\leq a_{k}\max\{\underset{t\geq 0}{\sup}y_{k-1}^{2}(t),K^{q_{k}}\}\\
        &\leq (a_{k})(a_{k-1})^{2}(a_{k-2})^{2^2}\cdots (a_1)^{2^{k-1}}\max\left\{\underset{t\geq 0}{\sup}y_{0}^{2^{k}}(t),\sum_{i=1}^{k-1}K^{q_{k-i}2^{i}}\right\}.
    \end{split}
\end{equation}
Since $K^{q_{k-i}2^{i}}<\widetilde{K}^{q_k}$, one concludes that
\begin{equation}
     y_{k}(t)\leq (2C\mu^{n+1})^{2^{k}-1}2^{(n+1)(2^{k+1}-k-2)}\max\left\{\underset{t\geq 0}{\sup}y_{0}^{2^{k}}(t),k\widetilde{K}^{q_k}\right\}.
\end{equation}
Taking both sides to the power of $q_{k}$ at the same time, we get
\begin{equation}\label{3.2.10}
\begin{split}
\Vert u \Vert_{L^{q_{k}}}&\leq \left(2C\mu^{n+1}\right)^{\frac{2^{k}-1}{q_{k}}}2^{\frac{(n+1)(2^{k+1}-k-2)}{q_{k}}}\max\left\{\underset{t\geq 0}{\sup}y_{0}^{\frac{2^{k}}{q_{k}}}(t),k^{\frac{1}{q_{k}}}\widetilde{K}\right\}\\
&\leq 2C\mu^{n+1}2^{2(n+1)}\max\{\underset{t\geq 0}{\sup}y_{0}(t),k^{\frac{1}{q_{k}}}\widetilde{K} \}.
\end{split}
\end{equation}
On the other hand, take $k=0$, and according to (\ref{3.1.1}), we know
\begin{equation}
    y_{0}(t)=\int_{\mathbb{R}^{n}}u^{q_0}dx=\int_{\mathbb{R}^{n}}u^{\mu}dx\leq C.
\end{equation}
Then the estimate is obtained by passing to the limit $k\to \infty$ in (\ref{3.2.10})
\begin{equation}
    \Vert u \Vert_{L^{\infty}(0,T;L^{\infty}(\mathbb{R}^{n}))}\leq C.
\end{equation}
\qed
\subsection{The case of \ \texorpdfstring{$0<b<1-\frac{\alpha}{n}$}{lg}}
\begin{lemma}\label{3.3.1.1}
Under the assumption of Theorem \ref{1.1.2}, the solution $(u, v)$ of system \ref{1.1} satisfies the following prior estimations
     \begin{flalign}\label{3.3.1}
        && \Vert u \Vert_{L^{\infty}(0,T;L^{\frac{n}{\alpha}}(\mathbb{R}^{n}))} &\leq C, \ \ \left\Vert D^{\frac{\alpha}{2}} u^{\frac{n}{2\alpha}} \right\Vert_{L^{2}(0,T;L^{2}(\mathbb{R}^{n}))}\leq C\\
       && \Vert \nabla v \Vert_{L^{\frac{n}{2}+1}(0,\infty;L^{r_1}(\mathbb{R}^{n}))} &\leq C,\  \frac{n}{n-1}<r_1\leq \frac{n(n+\alpha)}{\alpha(n-1)-n}.
       \end{flalign}
\end{lemma}
\begin{pro}
We also divide into two cases:
   \begin{enumerate} [(i)]
       \item $1<q\leq \frac{1}{1-b}$
       
       Similar to (i) of Lemma \ref{3.2}, we can see that (\ref{3.3.11}) is true.
       \item $q=\frac{n}{\alpha}>\frac{1}{1-b}$
        
       Similar to the discussion of Lemma \ref{3.1}, we can still get (\ref{3.1.3}) to hold. Applying the interpolation inequality and Sobolev inequality to $\int_{\mathbb{R}^{n}}u^{q+1}dx$, we obtain
       \begin{equation}\label{3.3.2}
           \begin{split}
               \Vert u \Vert_{L^{q+1}}^{q+1} &\leq \Vert u \Vert_{L^{\frac{(nq}{n-\alpha}}}^{(q+1)\tau_4}\Vert u \Vert_{L^{\frac{n}{\alpha}}}^{(q+1)(1-\tau_4)}\\
               &=\Vert u^{\frac{q}{2}} \Vert_{L^{\frac{2n}{n-\alpha}}}^{\frac{2\tau_4(q+1)}{q}}\Vert u \Vert_{L^{\frac{n}{\alpha}}}^{(q+1)(1-\tau_4)}\\
               &\leq S_{\alpha, n}^{\frac{2\tau_4(q+1)}{q}}\left\Vert D^{\frac{\alpha}{2}} u^{\frac{q}{2}} \right\Vert_{L^{2}}^{\frac{2(q+1)\tau_4}{q}}\Vert u \Vert_{L^{\frac{n}{\alpha}}}^{(q+1)(1-\tau_4)}.
          \end{split}
       \end{equation}
       where $\tau_4$ satisfies $\frac{1}{q+1}=\frac{\tau_4(n-\alpha)}{nq}+\frac{1-\tau_4}{\frac{n}{\alpha}}$, it is easy to see
       $\tau_4=\frac{\frac{1}{q+1}-\frac{\alpha}{n}}{\frac{n-\alpha}{nq}-\frac{\alpha}{n}}.$ 
       (\ref{3.3.2}) can be written as
       \begin{equation}\label{3.3.3}
           \Vert u \Vert_{L^{q+1}(\mathbb{R}^{n})}^{q+1}\leq S_{\alpha. n}^{2}\left\Vert D^{\frac{\alpha}{2}} u^{\frac{q}{2}}\right\Vert_{L^{2}}^{2}\Vert u \Vert_{L^{\frac{n}{\alpha}}(\mathbb{R}^{n})}.
       \end{equation}
       On the other hand, if $\sigma$ is small enough to satisfy $1<1+\sigma<\frac{1}{1-b}$, it can be known from the case of (i)
       \begin{equation}\label{3.3.4}
           \Vert u \Vert_{L^{\infty}(0,T;L^{{1+\sigma}}(\mathbb{R}^{n}))}\leq C.
       \end{equation}
       Next, interpolation inequality and Sobolev inequality are applied to $\Vert u \Vert_{L^{q}(\mathbb{R}^{n})}^{q}$
       \begin{equation}\label{3.3.5}
           \begin{split}
               \Vert u \Vert_{L^{q}(\mathbb{R}^{n})}^{q} &\leq \Vert u \Vert_{L^{\frac{nq}{n-\alpha}}}^{q(1-\tau_5)}\Vert u \Vert_{L^{1+\sigma}(\mathbb{R}^{n})}^{q\tau_5}\\
               &\leq S_{\alpha,n}^{2(1-\tau_5)}\left\Vert D^{\frac{\alpha}{2}} u^{\frac{q}{2}}\right\Vert_{L^{2}(\mathbb{R}^{n})}^{2(1-\tau_5)}\Vert u \Vert_{L^{1+\sigma}(\mathbb{R}^{n})}^{q\tau_5},
           \end{split}
       \end{equation}
       where $\tau_5$ satisfies $\frac{1}{q}=\frac{q(1-\tau_5)}{\frac{nq}{n-\alpha}}+\frac{\tau_5}{1+\sigma}$, calculated to get
       $$\tau_5=\frac{\frac{1}{q}-\frac{n-\alpha}{nq}}{\frac{1}{1+\sigma}-\frac{n-\alpha}{nq}}. $$
       applying Young's inequality to (\ref{3.3.5}) 
       \begin{equation}
           aq\Vert u \Vert_{L^{q}(\mathbb{R}^{n})}^{q}\leq \varepsilon_1\Vert D^{\frac{\alpha}{2}} u^{\frac{q}{2}} \Vert_{L^{2}(\mathbb{R}^{n})}^{2(1-\tau_5)l_1}+C(\varepsilon_1)(aq)^{l_2}S_{\alpha,n}^{2(1-\tau_5)l_2}\Vert u \Vert_{L^{1+\sigma}(\mathbb{R}^{n})}^{q\tau_5l_2}.
       \end{equation}
       Take $l_1=\frac{1}{1-\tau_5}$, $l_2=\frac{1}{1-\frac{1}{l_1}}=\frac{1}{\tau_5}$, then the following formula holds for any $\varepsilon_1 >0$
       \begin{equation}\label{3.3.6}
           aq\Vert u \Vert_{L^{q}(\mathbb{R}^{n})}^{q}\leq \varepsilon_1\Vert \nabla u \Vert_{L^{p}(\mathbb{R}^{n})}^{p}+C(\varepsilon_1)(aq)^{\frac{1}{\tau_5}}S_{\alpha,n}^{\frac{p(1-\tau_5)}{\tau_5}}\Vert u \Vert_{L^{1+\sigma}(\mathbb{R}^{n})}^{q}.
       \end{equation}
       Let $\varepsilon_1=\frac{2\nu(q+1)}{q}$, from (\ref{3.3.3}),(
       \ref{3.3.4}) and (\ref{3.3.6}), we get
       \begin{align}\label{3.3.7}
               &  \partial_{t}^{\beta}\int_{\mathbb{R}^{n}}u^{q}dx
               + \left(\frac{2(q-1)}{q}-S_{\alpha,n}^{2}(q-1-bq )\Vert u \Vert_{L^{\frac{n}{\alpha}}(\mathbb{R}^{n})}\right)\left\Vert D^{\frac{\alpha}{2}} u^{\frac{q}{2}}\right\Vert_{L^{2}(\mathbb{R}^{n})}^{2}\notag\\
               \leq{} & a^{\frac{1}{\tau_5}}C_{0},
       \end{align}
       where $C_0$ is a constant that depends on $q,S_{\alpha,n}, \Vert u_0 \Vert_{L^{1+\sigma}(\mathbb{R}^{n})}$. In addition, from the initial condition  $\Vert u_0 \Vert_{L^{\frac{n}{\alpha}}(\mathbb{R}^{n})}<\frac{C_{*}}{2}$, we obtain
       \begin{equation} \frac{2(q-1)}{q}-S_{\alpha,n}^{2}(q-1-bq)\Vert u_{0}\Vert_{L^{\frac{n}{\alpha}}(\mathbb{R}^{n})}
               > \frac{2(q-1)}{q}-S_{\alpha,n}^{2}(q-1-bq)\frac{C_{*}}{2} =:\sigma_1.
       \end{equation}
       Then take the small value $a$ satisfies $a<\left(\frac{(\frac{C_{*}}{2})^{\frac{n}{\alpha}}-\Vert u_0 \Vert_{L^{\frac{n}{\alpha}}(\mathbb{R}^{n})}^{\frac{n}{\alpha}}}{C_0T}\right)^{\tau_5}$, we obtain
       \begin{equation}
           \Vert u \Vert_{L^{\frac{n}{\alpha}}(\mathbb{R}^{n})} <\frac{C_{*}}{2},\ \ \ 0<t<T.
           \end{equation}
           That is, for any $0<t<T$
           $$\frac{2(q-1)}{q}S_{\alpha,n}^{2}(q-1-bq)\Vert u \Vert_{L^{\frac{n}{\alpha}(\mathbb{R}^{n})}} > \sigma_1.$$
           Then (\ref{3.3.7}) can be written as 
           \begin{equation}\label{3.3.8}
              \partial_{t}^{\beta}\int_{\mathbb{R}^{n}}u^{\frac{n}{\alpha}}dx+\sigma_{1}\Vert D^{\frac{\alpha}{2}}u^{\frac{n}{2\alpha}}\Vert_{L^{2}}^{2}\leq a^{\frac{1}{\tau_5}}C_0.
           \end{equation}
          By Lemma \ref{2.11}, we obtain
          \begin{equation}\label{3.3.9}
              \Vert u \Vert_{L^{\frac{n}{\alpha}}}^{\frac{n}{\alpha}}\leq \Vert u \Vert_{L^{q}}^{q}E_{\beta}(t^{\alpha})\leq C.
          \end{equation}
        From (\ref{3.3.8}) and (\ref{3.3.9}), it is easy to get (\ref{3.3.1}).
        Otherwise take $q=\frac{n}{\alpha}$ in (\ref{3.3.3}), combine (\ref{3.3.1}), we obtain
        \begin{equation}\label{3.3.9.1}
            \Vert u \Vert_{1+\frac{n}{\alpha}(0,T;L^{1+\frac{n}{\alpha}}(\mathbb{R}^{n}))}\leq C.
        \end{equation}
        
         The proof of $\nabla v$ can be obtained in a similar way to Lemma \ref{3.1}, and there is $\frac{n}{n-1}<r_1\leq \frac{n(n+\alpha)}{\alpha(n-1)-n}$ when $1<q\leq 1+\frac{n}{\alpha}$ is taken. That is, the conclusion is established.\qed
\end{enumerate}
\end{pro}
Next, we discuss the estimation of the spatial derivative and the estimation of the time derivative, which are two necessary conditions for applying the Aubin-Lions-Dubinski\u{\i} Lemma.
\begin{lemma}
Under the assumption of Theorem \ref{1.1.2}, the solution $(u, v)$ of system \ref{1.1} satisfies satisfies the following spatial derivative estimation
\begin{equation}\label{3.3.10}
    \left\Vert D^{\frac{\alpha}{2}}u\right\Vert_{L^{2}(0,T;L^{2}(\mathbb{R}^{n}))}\leq C.
\end{equation}
\end{lemma}
\begin{pro}
Multiplying the first equation in (\ref{1.1}) by $2u$ and integrating, we have
\begin{equation}\label{3.3.11}
                \int_{\mathbb{R}^{n}}2u\partial_{t}^{\beta}udx=-2\int_{\mathbb{R}^{n}}\left\lvert D^{\frac{\alpha}{2}}u\right\rvert^{2}dx+(1-2b)\int_{\mathbb{R}^{n}}u^{3}dx+2a\int_{\mathbb{R}^{n}}u^{2}dx.
            \end{equation}
\end{pro}
(1) $1-2b\leq 0$\\
(\ref{3.3.11}) can be written as 
                \begin{equation}
             \int_{\mathbb{R}^{n}}2u\partial_{t}^{\beta}udx\leq -2\Vert D^{\frac{\alpha}{2}}u \Vert_{L^{2}(\mathbb{R}^{n})}+a\Vert u \Vert_{L^{2}(\mathbb{R}^{n})}.
                \end{equation}
                It is easy to see that
                \begin{equation}
                    _{0}^{c}\textrm{D}_{t}^{\alpha}\Vert u \Vert_{L^{2}(\mathbb{R}^{n})}^{2}+\Vert D^{\frac{\alpha}{2}}u \Vert_{L^{2}(\mathbb{R}^{n})}^{2}\leq a \Vert u \Vert_{L^{2}(\mathbb{R}^{n})}^{2}.
                \end{equation}
                Then (\ref{3.3.10}) can be obtained from Lemma \ref{3.3.1.1}.\\
(2) $1-2b>0$ and $n=2\alpha$
By (\ref{3.3.1}) and (\ref{3.3.4}), we can get $\int_{\mathbb{R}^{n}}u^2dx\leq C$, in addition, take $q=2$ in (\ref{3.3.3}), and get
                \begin{equation}\label{3.3.12}
                    \Vert u \Vert_{L^{3}(\mathbb{R}^{n})}\leq S_{\alpha,n}^{2}\left\Vert D^{\frac{\alpha}{2}} u \right\Vert_{L^{2}(\mathbb{R}^{n})}^{2}\Vert u \Vert_{L^{\frac{n}{\alpha}}(\mathbb{R}^{n})}\leq C.
                \end{equation}
                Substituting (\ref{3.3.12}) into (\ref{3.3.11}) yields to
                \begin{equation}\label{3.3.13}
                    \partial_{t}^{\alpha}\Vert u \Vert_{L^{2}(\mathbb{R}^{n})}^{2}\leq -\left(2-(1-2b)S_{\alpha,n}^{2}\Vert u \Vert_{L^{\frac{n}{\alpha}}(\mathbb{R}^{n})}\right)\left\Vert D^{\frac{\alpha}{2}} u \right\Vert_{L^{2}(\mathbb{R}^{n})}^{2}+C.
                \end{equation}
                From $\Vert u \Vert_{L^{\frac{n}{\alpha}(\mathbb{R}^{n})}}<\frac{C_{*}}{2}$, we know $2-(1-2b)S_{\alpha,n}^{2}\Vert u \Vert_{L^{\frac{n}{\alpha}}(\mathbb{R}^{n})}>0$, i.e.,
                \begin{equation}                    \left(2-(1-2b)S_{\alpha,n}^{2}\Vert u\Vert_{L^{\frac{n}{\alpha}}(\mathbb{R}^{n})}\right)\left\Vert D^{\frac{\alpha}{2}} u \right\Vert_{L^{2}(\mathbb{R}^{n})}^{2}\leq C.
                \end{equation}
                That is, (\ref{3.3.10}) holds.\\
(3) $1-2b>0$ and $n> 2\alpha$
By (\ref{3.3.1}) and (\ref{3.3.4}), we have
                $$\Vert u \Vert_{L^{2}(\mathbb{R}^{n})}^{2}\leq C,\ \ \ \Vert u \Vert_{L^{3}(\mathbb{R}^{n})}^{3}\leq C.$$
                That is, we have
                \begin{equation}
                    \partial_{t}^{\alpha}\int_{\mathbb{R}^{n}}u^{2}dx+2\int_{\mathbb{R}^{n}}\left\lvert D^{\frac{\alpha}{2}}\right\rvert^{2}\leq C.
                \end{equation}
                (\ref{3.3.10}) is easy to get from the above formula.
                
                Therefore, the above proof shows that (\ref{3.3.10}) holds.
\begin{lemma}
Under the assumption of Theorem \ref{1.1.2}, the solution $(u, v)$ of system \ref{1.1} satisfies satisfies the following time derivative estimation
\begin{equation}\label{3.3.14}
     \Vert \partial_{t}^{\beta}u\Vert_{L^{2}(0,T;W^{-\alpha,r_2}(U))}\leq C,\ \ r_{2}=\min\{2,\frac{n(n+\alpha)}{n\alpha+n+\alpha}\}.
\end{equation}
\end{lemma}
\begin{pro}
 We first prove that the following estimator holds:
       \begin{equation}\label{3.3.15}
       \Vert u\nabla v \Vert_{L^{1+\frac{n}{\alpha}}(0,T;L^{\frac{n(n+\alpha)}{n\alpha+n+\alpha}}(\mathbb{R}^{n}))}\leq C.
   \end{equation}
By H$\ddot{o}$lder's inequality
    \begin{equation}
        \int_{\mathbb{R}^{n}}u^{\frac{n(n+\alpha)}{n\alpha+n+\alpha}}\lvert \nabla v \rvert^{\frac{n(n+\alpha)}{n\alpha+n+\alpha}}dx\leq \left(\int_{\mathbb{R}^{n}}u^{\frac{n(n+\alpha)}{n\alpha+n+\alpha}l_5}dx\right)^{\frac{1}{l_5}}\left(\int_{\mathbb{R}^{n}}\lvert \nabla v \rvert^{\frac{n(n+\alpha)}{n\alpha+n+\alpha}l_6}dx\right)^{\frac{1}{l_6}}
    \end{equation}
    where $\frac{1}{l_5}+\frac{1}{l_6}=1$, take $\frac{n(n+\alpha)}{n\alpha+n+\alpha}l_6=n$,we get $l_6=\frac{n\alpha+n+\alpha}{n+\alpha}$,so we have $l_5=\frac{n\alpha+n+\alpha}{n\alpha}$, we obtain
    \begin{equation}
        \Vert u \nabla v\Vert_{L^{\frac{n(n+\alpha)}{n\alpha+n+\alpha}}(\mathbb{R}^{n})}\leq \Vert u \Vert_{L^{\frac{n}{\alpha}+1}(\mathbb{R}^{n})}\Vert \nabla v\Vert_{L^{n}(\mathbb{R}^{n})}
    \end{equation}
From (\ref{3.3.9.1}) and (\ref{3.3.1}), (\ref{3.3.15}) holds. 

Take test function $\varphi(x,t)$ with $p_{1}>1$ and
$$\Vert \varphi \Vert_{L^{q_1}\left((0,T);W^{\alpha,r_{2}^{*}}(\mathbb{R}^{n})\right)},$$
where $r_{2}^{*}$ is the conjugate index of $r_2$, $r_2=\min\left\{\frac{n}{\alpha},\frac{n(n+\alpha)}{n\alpha+n+\alpha}\right\}$. we obtain
\begin{equation}
\begin{split}
    \int_{\mathbb{R}^{n}}\partial_{t}^{\beta}u\phi dx&=- \int_{\mathbb{R}^{n}} uD^{\alpha}\phi dx+ \int_{\mathbb{R}^{n}}u\nabla \cdot \nabla \phi dx\\
    &+a \int_{\mathbb{R}^{n}}u\phi dx-b \int_{\mathbb{R}^{n}}u^{2}\phi dx.
\end{split}
\end{equation}
Using H$\ddot{o}$lder's inequality
\begin{equation}
    \begin{split}
    \langle  \partial_{t}^{\beta} u, \phi \rangle &\leq \langle u\nabla v, \nabla \varphi \rangle +\langle u, D^{\alpha}\varphi \rangle+a\langle u,\varphi \rangle +b\langle u^2,\varphi\rangle\\
    &\leq  \Vert  u\Vert_{L^{\frac{n}{\alpha}}(\mathbb{R}^{n})}\Vert D^{\alpha}\phi \Vert_{L^{\frac{n}{n-\alpha}}(\mathbb{R}^{n})}+\Vert u\nabla v\Vert_{L^{\frac{n(n+\alpha)}{n\alpha+n+\alpha}}(\mathbb{R}^{n})}\Vert \nabla \phi \Vert_{L^{\frac{n(n+\alpha)}{n^{2}-n+\alpha}}(\mathbb{R}^{n})}\\
    &+ a\Vert u \Vert_{L^{\frac{n}{\alpha}}(\mathbb{R}^{n})}\Vert \phi \Vert_{L^{\frac{n}{n-\alpha}}(\mathbb{R}^{n})}+b\Vert u^{2}\Vert_{L^{\frac{n+\alpha}{2\alpha}}(\mathbb{R}^{n})}\Vert \phi\Vert_{L^{\frac{n+2}{n-2}}(\mathbb{R}^{n})}\Vert \phi \Vert_{L^{\frac{n+\alpha}{n-\alpha}}(\mathbb{R}^{n})}.
    \end{split}
\end{equation}
 Thus for any $T>0$, when $t\in[0,T]$, we have
\begin{equation}
    \begin{split}
       \int_{0}^{T} \Vert \partial_{t}^{\beta}u\Vert_{W^{-\alpha,r_2}(\mathbb{R}^{n})} 
        dt&\leq C\left({\begin{array}{c}
        \int_{0}^{T}\Vert u \Vert_{L^{\frac{n}{\alpha}}(\mathbb{R}^{n})}dt
        +\int_{0}^{T}\Vert u\nabla v\Vert_{L^{\frac{n(n+\alpha)}{n\alpha+n+\alpha}}(\mathbb{R}^{n})}^{2}dt\\
        + \int_{0}^{T}\Vert  u\Vert_{L^{\frac{n}{\alpha}}(\mathbb{R}^{n})}dt+ \int_{0}^{T}\Vert u \Vert_{L^{1+\frac{n}{\alpha}}(\mathbb{R}^{n})}dt\\
        \end{array}}\right).\\
        &\leq 3C(T)+C\int_{0}^{T}\Vert u \Vert_{L^{1+\frac{n}{\alpha}}}\Vert \nabla v\Vert_{L^{n}}dt\\
        &\leq 4C(T).
    \end{split}
    \end{equation}
From (\ref{3.3.15}), that is 
$$\Vert \partial_{t}^{\beta}u\Vert_{L^{q_1}(0,T;W^{-\alpha,r_2}(U))}\leq C,\ \ r_2=\min\{\frac{n}{\alpha},\frac{n(n+\alpha)}{n\alpha+n+\alpha}\},$$
where $q_{1}>1$ and $\frac{1}{q_1}+\frac{1}{p_1}=1$. \qed
\end{pro}

\noindent $\mathbf{Proof\ of\ Theorem\ref{1.1.2}:}$\\
 \textbf{Step 1.}(The regularization problem and a priori estimates) Similar to the previous discussion, the solution of the regularization equation (\ref{3.1.5}) satisfies all the estimators of the above lemmas.\\
  \begin{equation}
  \Vert u_{\varepsilon} \Vert_{L^{\infty}(0,T;L^{1}(\mathbb{R}^{n})\cap L^{\frac{n}{\alpha}}(\mathbb{R}^{n}))} \leq C,
    \end{equation}
    \begin{equation}
    \left\Vert D^{\frac{\alpha}{2}}u_{\varepsilon}\right\Vert_{L^{2}(0,T;L^{2}(\mathbb{R}^{n}))}\leq C.
    \end{equation}
    \begin{equation}
     \Vert \partial_{t}^{\beta}u_{\varepsilon}\Vert_{L^{2}(0,T;W^{-\alpha,r_2}(U))}\leq C,\ \ r_{2}=\min\{\frac{n}{\alpha},\frac{n(n+\alpha)}{n\alpha+n+\alpha}\},
    \end{equation}
    \begin{equation}
       \Vert \nabla v_{\varepsilon} \Vert_{L^{\frac{n}{\alpha}+1}(0,\infty;L^{r_1}(\mathbb{R}^{n}))} \leq C, \frac{n}{n-1}<r_1\leq \frac{n(n+\alpha)}{\alpha(n-1)-n}.
    \end{equation}
 \textbf{Step 2.}(Compactness)
 Similar to the discussion of Theorem \ref{1.1.1}, we also have
  \begin{equation}
  \Vert u_{\varepsilon} \Vert_{L^{q}(0,T;P_{+}(\Omega)} \leq C,
    \end{equation}
    \begin{equation}
    \left\Vert u_{\varepsilon}\right\Vert_{L^{q}(0,T;L^{q}(\Omega))}\leq C.
    \end{equation}
    \begin{equation}
     \Vert \partial_{t}^{\beta}u_{\varepsilon}\Vert_{L^{2}(0,T;W^{-\alpha,r_2}(\Omega))}\leq C,\ \ r_{2}=\min\{\frac{n}{\alpha},\frac{n(n+\alpha)}{n\alpha+n+\alpha}\}.
    \end{equation}
    Using the Aubin-Lions-Dubinski\u{\i} Lemma, we deduce that there are subsequences of $\{u_{\varepsilon}\}, \{v_{\varepsilon}\}$ (still denoted as $\{u_{\varepsilon}\}, \{v_{\varepsilon}\}$), and $u, v$ satisfying the regularization property
    \begin{gather}
        u \in L^{\infty}(0,T;L^{1}(\mathbb{R}^{n}) \cap L^{\frac{n}{\alpha}}(\mathbb{R}^{n}))\cap L^{2}(0,T;H^{\frac{\alpha}{2}}(\mathbb{R}^{n})),\\
        D^{\frac{\alpha}{2}} u \in L^{2}(0,T;L^2(\mathbb{R}^{n})),\\
        \partial_{t}^{\beta}u \in L^{q_{1}}(0,T;W^{-\alpha,r_2}(\mathbb{R}^{n})),\ \ r_2 = \min\{\frac{n}{\alpha},\frac{n(n+\alpha)}{n\alpha+n+\alpha}\},\\
        \nabla v \in L^{\frac{n}{\alpha}+1}((0,\infty);L^{s_1}(\mathbb{R}^{n})),\ \ \frac{n}{n-1}<r_1\leq \frac{n(n+\alpha)}{\alpha(n-1)-n}.
     \end{gather}
    so that the following weak convergence relation holds
     \begin{gather}
        u_{\varepsilon}\rightharpoonup u\ \ {\rm weakly^{*}}\  {\rm in} \  L^{\infty}(0,T;L^{1}(\mathbb{R}^{n})\cap L^{\frac{n}{\alpha}}(\mathbb{R}^{n}))\cap L^{2}(0,T;H^{\frac{\alpha}{2}}(\mathbb{R}^{n})),\\
         D^{\frac{\alpha}{2}} u_{\varepsilon} \rightharpoonup   D^{\frac{\alpha}{2}}\ \  {\rm weakly}\  {\rm in}\  L^{2}(0,T;L^{2}(\mathbb{R}^{n})),\\
         \partial_{t}^{\beta}u_{\varepsilon}\rightharpoonup  \partial_{t}^{\beta}u \ \  {\rm weakly}\ {\rm in} \ L^{q_1}(0,T;W^{-\alpha,r_2}(\mathbb{R}^{n})), r_2 = \min\{\frac{n}{\alpha},\frac{n(n+\alpha)}{n\alpha+n+\alpha}\},\\
        \nabla v_{\varepsilon}\rightharpoonup \nabla v \ \ {\rm weakly^{*}}\ {\rm in}\  L^{\frac{n}{\alpha}+1}(0,T;L^{s_2}(\mathbb{R}^{n})),\ \ \ \ \frac{n}{n-1}<r_1\leq \frac{n(n+\alpha)}{\alpha(n-1)-n}.
    \end{gather}
The same can be obtained
 \begin{equation}
     u_{\varepsilon}\to u {\ \rm\ in\ } L^{q}(0,T;L^{q}(\Omega)\ \ {\rm\ as\ }\varepsilon\to 0.
 \end{equation}
 The following consistent strong convergence holds
 \begin{equation}
     u^{\varepsilon}\to u \ {\ \rm\ in\ } L^{q}(0,T;L^{q}(B_{k})\ \ {\rm\ as\ }\varepsilon\to 0,\ \forall k.
 \end{equation}
    \textbf{Step 3.} (The existence of global weak solutions)Then for $\forall \phi \in C_{c}^{\infty}(\mathbb{R}^{n})$, we have
     \begin{align}
        &\int_{0}^{T}\int_{\mathbb{R}^{n}}(u_{\varepsilon}-u_0)\phi(x,t)\tilde{\partial}_{T}^{\beta}\phi(x,t)dxdt\notag\\
        ={}& - \int_{0}^{T}\int_{\mathbb{R}^{n}}\left[u_{\varepsilon}(x,t)D^{\alpha}\phi(x)\right]dxdt+\int_{0}^{T}\int_{\mathbb{R}^{n}}u_{\varepsilon}(x) \nabla v_{\varepsilon}(x)\cdot \nabla \phi(x)dxdt\notag\\
        +{} & a\int_{0}^{T}\int_{\mathbb{R}^{n}}u_{\varepsilon}(x)\phi(x)dxdt-b\int_{0}^{T}\int_{\mathbb{R}^{n}}u_{\varepsilon}^{2}(x)\phi(x)dxdt.
\end{align}
    Let $\varepsilon \to 0$, Using the above convergence relationship, the remaining proof steps are the same as the proof process of Theorem \ref{1.1.1} (i). The limit functions $u,v$ can be obtained as a weak solution of equation (\ref{1.1}).\qed
\section{Uniqueness of weak solutions}
First we pick $A=(-\Delta)^{\frac{\alpha}{2}}$($1<\alpha\leq2$), and consider operators $T_{\alpha}^{\beta}(t)$ defined by
$$f(x)\mapsto T_{\alpha}^{\beta}(t)f(x):=t^{\beta-1}E_{\beta,\beta}(-t^{\beta}A)f(x).$$
\begin{lemma}\cite{li2018cauchy}\label{4.1}
Let $0<\beta<1$ and $1<\alpha\leq2$. Then\\
(1) Let $r\in [1,\infty)$. We define $\zeta_1=\frac{nr}{n-2r\alpha}$ if $n>2r\alpha$ and $\zeta_1=\infty$ otherwise. Then, for any $p\in [1,\zeta_1)$, we have
\begin{equation}
\left\Vert T_{\alpha}^{\beta}(t)u\right\Vert_{L^{p}}\leq Ct^{-\frac{n\beta}{\alpha}(\frac{1}{r}-\frac{1}{p})+\beta-1}\Vert u\Vert_{L^{r}}.
\end{equation}
(2) Let $r\in[1,\infty)$, define $\zeta_2=\frac{nr}{n+r(1-2\alpha)}$ if $n>r(2\alpha-1)$, and $\zeta_2=\infty$ otherwise. Then for $p=[r,\zeta_2)$ there is $C>0$ satisfying
\begin{equation}
    \left\Vert \nabla T_{\alpha}^{\beta}(t)u\right\Vert_{p}\leq Ct^{-\frac{n\beta}{\alpha}(\frac{1}{r}-\frac{1}{p})-\frac{\beta}{\alpha}+\beta-1}\Vert u\Vert_{r}
\end{equation}
\end{lemma}
\begin{lemma}\label{4.1.1.1}
(i) Let $u$ be the weak solution of model (\ref{1.1}) under the condition of Theorem \ref{1.1.1}. For any $q\geq \frac{n}{\alpha}$, there exists a constant $C(a, q, n)$ such that the following equation holds
\begin{equation}\label{4.4}
    t^{q-\frac{n}{\alpha}}\Vert u \Vert_{q}^{L^{q}(\mathbb{R}^{n})}\leq C(a, q, n),\ \ \forall 0<t<T,
\end{equation}
where $C(a, q, n)$ represents a constant that depends only on $a, q, n$.\\
(ii) Suppose further that for any given $0<\varepsilon<1$, $\Vert u_{0}\Vert_{L^{\frac{n}{\alpha}+\varepsilon}}<\infty$. Then for any $q\geq \frac{n}{\alpha}+\varepsilon$ we have
\begin{equation}\label{4.5}
     t^{q-\frac{n}{\alpha}-\frac{\alpha q\varepsilon}{n+\alpha\varepsilon}}\Vert u \Vert_{q}^{L^{q}(\mathbb{R}^{n})}\leq C(a, q, n),\ \ \forall 0<t<T.
\end{equation}
\end{lemma}
\begin{pro}
From (\ref{3.1.4}) and (\ref{3.2.6.1}), we know that under the condition of satisfying Theorem \ref{1.1.1}, the  equations is finally written as
\begin{equation}
    \partial_{t}^{\alpha} \int_{\mathbb{R}^{n}}u^{q}dx+C_7\int_{\mathbb{R}^{n}}\left\lvert D^{\frac{\alpha}{2}}u^{\frac{q}{2}}\right\rvert^{2}dx\leq aq\int_{\mathbb{R}^{n}}u^{q}dx+C.
\end{equation}
Combine (\ref{3.3.4}), (\ref{3.3.5}) and (\ref{3.3.6}), we obtain
\begin{equation}\label{4.2}
    \partial_{t}^{\alpha} \int_{\mathbb{R}^{n}}u^{q}dx+C_8\int_{\mathbb{R}^{n}}\left\lvert D^{\frac{\alpha}{2}}u^{\frac{q}{2}}\right\rvert^{2}dx\leq C,
\end{equation}
where $C_8$ is a constant. Next using Interpolation inequality to $\Vert u \Vert_{L^{q}}$
\begin{equation}\label{4.2.1}
    \Vert u \Vert_{L^{q}}\leq \Vert u \Vert_{L^{\frac{n}{\alpha}}}^{\tau_6}\Vert u \Vert_{L^{\frac{nq}{n-\alpha}}}^{1-\tau_6},
\end{equation}
where $0<\tau_6<1$ and satisfies $\frac{1}{q}=\frac{\alpha\tau_6}{n}+\frac{(1-\tau_6)(n-\alpha)}{nq}$, i.e. $\tau_6=\frac{\alpha}{\alpha q-n+\alpha}$.
we have
\begin{equation}\label{4.3}
\begin{split}
    \Vert u \Vert_{L^{q}}&\leq \Vert u \Vert_{L^{\frac{n}{\alpha}}}^{\frac{\alpha}{\alpha q-n+\alpha}}\Vert u \Vert_{L^{\frac{nq}{n-\alpha}}}^{\frac{\alpha q-n\alpha}{\alpha q-n+\alpha}}\\
    &\leq \Vert u \Vert_{L^{\frac{n}{\alpha}}}^{\frac{\alpha}{\alpha q-n+\alpha}}S_{\alpha,n}^{\frac{2\alpha(q-n)}{\alpha q-n+\alpha}}\left \Vert D^{\frac{\alpha}{2}}u^{\frac{q}{2}}\right\Vert_{L^2}^{\frac{2\alpha(q-n)}{\alpha q-n+\alpha}}\\
    &=\bar{C}\left \Vert D^{\frac{\alpha}{2}}u^{\frac{q}{2}}\right\Vert_{L^2}^{\frac{2\alpha(q-n)}{\alpha q-n+\alpha}},
\end{split}
\end{equation}
where $\bar{C}=S_{\alpha,n}^{\frac{2\alpha(q-n)}{\alpha q-n+\alpha}}\Vert u \Vert_{L^{\frac{n}{\alpha}}}^{\frac{\alpha}{\alpha q-n+\alpha}}$. Substitute (\ref{4.3}) into (\ref{4.2})
\begin{equation}
    _{0}^{c}\textrm{D}_{t}^{\beta}\int_{\mathbb{R}^{n}}u^{q}dx+\frac{C_8}{\bar{C}}\Vert u \Vert_{L^{q}}^{\frac{q(\alpha q-n+\alpha)}{\alpha (q-n)}}\leq C.
\end{equation}
By Lemma \ref{2.12}, it is easy to get (\ref{4.4}).

Take $q=\frac{n}{\alpha}+\varepsilon$ in (\ref{4.2}) to get
$\Vert u \Vert_{L^{\frac{n}{\alpha}+\varepsilon}}\leq C.$ Replace $\Vert u \Vert_{L^{\frac{n}{\alpha}}}$ with $\Vert u \Vert_{L^{\frac{n}{\alpha}+\varepsilon}}$ for interpolation in (\ref{4.2.1}), and (\ref{4.5}) is established in the same way.\qed
\end{pro}

\noindent $\mathbf{Proof\ of\ Theorem\ref{1.1.3}}$\\
 Suppose we assume that $(u_1,v_1)$ and $(u_2,v_2)$ are any two solutions of given system eq. (\ref{1.1}). Take $\tilde{u}=u_1-u_2$, $\tilde{v}=v_1-v_2$, then $(\tilde{u},\tilde{v})$ satisfy,
    \begin{equation}
        \left\{
        \begin{aligned}
    \partial_{t}^{\beta}\tilde{u} & = -(-\Delta)^{\frac{\alpha}{2}} \tilde{u} - \nabla \cdot (\tilde{u}\nabla v_1) - \nabla \cdot (u_2\nabla \tilde{v}) +  a\tilde{u} - bu_{1}^2 + bu_{2}^2,  \\
    -\Delta \tilde{v} & = \tilde{u},\\
    \tilde{u}(x,0) & = 0.\\
    \end{aligned}
        \right.
    \end{equation}
    Set
    $$M_{k,r}(t) = \underset{0<s<t}{\sup}s^{k}\Vert \tilde{u} \Vert_{L^{r}(\mathbb{R}^{n})}$$
    $$Q_{i,n}(t) = \underset{0<s<t}{\sup}s^{1-\frac{n\beta}{\alpha q}} \Vert u_{i} \Vert_{L^{q}(\mathbb{R}^{n})},\ \ i=1, 2,$$
    where $r > 1, q > \frac{n}{\alpha},  1<k <1+ \frac{n\beta}{\alpha q}$.
    Formally taking the Laplacian transform of (\ref{1.1}), we find that $u$ satisfies the following Duhamel type integral equation
    \begin{equation}
        \begin{split}
            \tilde{u}&=
            E_{\beta}(-t^{\beta}A)\tilde{u_0}-\beta\int_{0}^{t}(t-s)^{\beta-1}E_{\beta}^{'}(-(t-s)^{\beta}A)\left[{\begin{array}{cc}
            \nabla \cdot (\tilde{u}\nabla v_1)+\nabla \cdot (u_2\nabla \tilde{v})\\
            -a\tilde{u}+bu_{1}^{2}-bu_{2}^{2}\\
            \end{array}}\right]ds\\
            &=-\int_{0}^{t}(t-s)^{\beta-1}E_{\beta,\beta}(-(t-s)^{\beta}A)\left[{\begin{array}{cc}
            \nabla \cdot (\tilde{u}\nabla v_1)+\nabla \cdot (u_2\nabla \tilde{v})\\
            -a\tilde{u}+bu_{1}^{2}-bu_{2}^{2}\\
            \end{array}}\right]ds\\
            &=-\int_{0}^{t}(t-s)^{\beta-1}E_{\beta,\beta}(-(t-s)^{\beta}A)(\nabla \cdot (\tilde{u}\nabla v_1))ds\\
            &-\int_{0}^{t}(t-s)^{\beta-1}E_{\beta,\beta}(-(t-s)^{\beta}A)(\nabla \cdot (u_2\nabla \tilde{v}))ds\\
            &+a\int_{0}^{t}(t-s)^{\beta-1}E_{\beta,\beta}(-(t-s)^{\beta}A)\tilde{u}ds\\
            &+b\int_{0}^{t}(t-s)^{\beta-1}E_{\beta,\beta}(-(t-s)^{\beta}A)(u_{1}^{2}-u_{2}^{2})ds\\
            &\overset{\Delta}{=} I_1+I_2+I_3+I_4.
        \end{split}
    \end{equation}
    Next we estimate $I_1,I_2,I_3,I_4$. First  estimate $t^{k}\Vert I_1 \Vert_{L^{r}(\mathbb{R}^{n})}$, we obtain
    \begin{equation}
    \begin{split}
        t^{k}\Vert I_1\Vert_{L^{r}(\mathbb{R}^{n})} &=t^{k}\Vert \int_{0}^{t}T_{\alpha}^{\beta}(t-s)(\nabla \cdot (\tilde{u}\nabla v_{1}))ds\Vert_{L^{r}(\mathbb{R}^{n})}\\
        &\leq Ct^{k}\int_{0}^{t}\left\Vert \nabla T_{\alpha}^{\beta} (t-s) \tilde{u}\nabla v_{1}\right\Vert_{L^{r}(\mathbb{R}^{n})}ds\\
        & \leq Ct^{k}\int_{0}^{t}(t-s)^{-\frac{n\beta}{\alpha}(\frac{1}{\sigma}-\frac{1}{r})+\beta-1}\Vert \tilde{u}\nabla v_{1}\Vert_{L^{\sigma}(\mathbb{R}^{n})}ds\\
        & \leq Ct^{k}\int_{0}^{t}(t-s)^{-\frac{n\beta}{\alpha}(\frac{1}{\sigma}-\frac{1}{r})+\beta-1}\Vert u \Vert_{L^{r}(\mathbb{R}^{n})}\Vert \nabla v\Vert_{L^{r'}(\mathbb{R}^{n})}ds,
    \end{split}
    \end{equation}
where $r'>1$ satisfies $\frac{1}{\sigma}=\frac{1}{r}+\frac{1}{r'}$. Applying the weak Young's inequality to $\nabla v_{1}$
\begin{equation}\label{11111}
    \begin{split}
        \Vert \nabla v_1 \Vert_{L^{r'}(\mathbb{R}^{n})}&=\Vert F* u \Vert_{L^{r'}(\mathbb{R}^{n})}\\
        &\leq C\left\Vert \frac{x}{\lvert x \rvert^{n}}* u \right\Vert_{L^{r'}(\mathbb{R}^{n})}\\
        &\leq C\left\Vert \frac{x}{\lvert x \rvert^{n}}\right\Vert_{L_{\omega}^{p}(\mathbb{R}^{n})}\Vert u_{1}\Vert_{L^{q}(\mathbb{R}^{n})}\\
        &\leq C\Vert u_{1}\Vert_{L^{q}(\mathbb{R}^{n})},
    \end{split}
\end{equation}
where $p=\frac{n\beta}{n\beta-\alpha\beta-\alpha}$, $q>\frac{n}{\alpha}$ satisfy $1+\frac{1}{r'}=\frac{1}{p}+\frac{1}{q}$, so we have
\begin{align}
\begin{split}
    & t^{k}\Vert I_{1}\Vert_{L^{r}(\mathbb{R}^{n})}\notag\\
    \leq{} & Ct^{k}\int_{0}^{t}(t-s)^{-\frac{n\beta}{\alpha}(\frac{1}{\sigma}-\frac{1}{r})+\beta-1}\Vert \tilde{u}\Vert_{L^{r}(\mathbb{R}^{n})}\Vert u_{1} \Vert_{L^{q}(\mathbb{R}^{n})}ds\notag\\
    \leq{}&  Ct^{k}\int_{0}^{t}(t-s)^{-\frac{n\beta}{\alpha q}}s^{k}\Vert\tilde{u}\Vert_{L^{r}(\mathbb{R}^{n})}\int_{0}^{t}s^{1-\frac{n\beta}{\alpha q}}\Vert u_{1} \Vert_{L^{q}(\mathbb{R}^{n})}s^{\frac{n\beta}{\alpha q}-1-k}ds\notag\\
    \leq{}& CM_{k,q}(t)Q_{1,n}(t)\int_{0}^{1}\left(1-\frac{s}{t}\right)^{-\frac{n\beta}{\alpha q}}\left(1-\frac{s}{t}\right)^{{\frac{n\beta}{\alpha q}-1-k}}d\left(\frac{s}{t}\right).
\end{split}
\end{align}
Let $\rho =\frac{s}{t}$, then the above formula can be written as
\begin{equation}
    t^{k}\Vert I_{1}\Vert_{L^{r}(\mathbb{R}^{n})}\leq CM_{k,q}(t)Q_{1,n}(t)\int_{0}^{1}(1-\rho)^{\frac{\beta}{\alpha}-\frac{n\beta}{\alpha q}-1}\rho^{\frac{n\beta}{\alpha q}-\frac{\beta}{\alpha}-k}d\rho,
\end{equation}
where $q>\frac{n}{\alpha}$ and $1<k<1+\frac{n\beta}{\alpha q}$, so there is $\varpi \overset{\Delta}{=}1-\frac{n\beta}{\alpha q}>0$,$\varrho\overset{\Delta}{=}\frac{n\beta}{\alpha q}-k>0$, it is easy to see that
\begin{equation}
   \int_{0}^{1}(1-\rho)^{\frac{\beta}{\alpha}-\frac{n\beta}{\alpha q}-1}\rho^{\frac{n\beta}{\alpha q}-\frac{\beta}{\alpha}-k}d\rho=\mathcal{B}(\varpi,\varrho) <\infty.
\end{equation}
Hence there is
\begin{equation}\label{4.1.1}
    t^{k}\Vert I_1 \Vert_{L^{r}(\mathbb{R}^{n})}\leq C M_{k,q}(t)Q_{1,n}(t).
\end{equation}

We use a similar method to estimate $t^{k}\Vert I_2 \Vert_{L^{r}(\mathbb{R}^{n})}$.
\begin{equation}
\begin{split}
    t^{k}\Vert I_2\Vert_{L^{r}}&=t^{k}\left\Vert \int_{0}^{t}T_{\alpha}^{\beta}(t-s)\left(\nabla \cdot (u_2\nabla \tilde{v}))ds\right) \right\Vert_{L^{r}}\\
    &\leq Ct^{k}\int_{0}^{t}\left \Vert \nabla T_{\alpha}^{\beta}(t-s)(u_2\nabla \tilde{v})\right\Vert_{L^{r}}ds\\
    &\leq Ct^{k}\int_{0}^{t}(t-s)^{-\frac{n\beta}{\alpha}(\frac{1}{\sigma}-\frac{1}{r})+\beta-1}\Vert u_2\nabla \tilde{v}\Vert_{L^{\sigma}}ds\\
    &\leq Ct^{k}\int_{0}^{t}(t-s)^{-\frac{n\beta}{\alpha}(\frac{1}{\sigma}-\frac{1}{r})+\beta-1}\Vert u_2\Vert_{L^{q}}\Vert \nabla \tilde{v}\Vert_{L^{r'}}ds,
\end{split}
\end{equation}
where $1\leq \sigma \leq r\leq \infty$, $q>\frac{n}{\alpha}$ and $\forall r'>1$ satisfy $\frac{1}{\sigma}=\frac{1}{q}+\frac{1}{r'}$. Apply weak Young's inequality to $\nabla \tilde{v}$ can be obtained
\begin{equation}
        \begin{split}
            \Vert \nabla \tilde{v} \Vert_{L^{r^{'}}(\mathbb{R}^{n})}
            &\leq C \left\Vert \frac{x}{\lvert x \rvert^{n}} * (u_1-u_2) \right\Vert_{L^{r^{'}}(\mathbb{R}^{n})}\\
            &\leq C \left\Vert \frac{x}{\lvert x \rvert^{n}} \right\Vert_{L_{\omega}^{p}(\mathbb{R}^{n})} \Vert \tilde{u} \Vert_{L^{r}(\mathbb{R}^{n})}\\
            &\leq C\Vert \tilde{u} \Vert_{L^{r}(\mathbb{R}^{n})},
        \end{split}
    \end{equation}
where $\frac{1}{p}=\frac{n\beta-\alpha\beta-\alpha}{n\beta}$ satisfies $1+\frac{1}{r'}=\frac{1}{p}+\frac{1}{r}$
\begin{align}
    & t^{k}\Vert I_2 \Vert_{L^{r}}\notag\\
    \leq{} & Ct^{k}\int_{0}^{t}(t-s)^{-\frac{n\beta}{\alpha q}}\Vert u_2\vert_{L^{q}}\Vert \tilde{u}\Vert_{L^{r}}ds\notag\\
    \leq{} & Ct^{k}\int_{0}^{t}t^{-\frac{n\beta}{\alpha q}}\left(1-\frac{s}{t}\right)^{-\frac{n\beta}{\alpha q}}s^{k}\Vert \tilde{u}\Vert_{L^{r}}s^{1-\frac{n\beta}{\alpha q}}\Vert u_{1}\Vert_{L^{q}}s^{\frac{n\beta}{\alpha q}-1-k}ds\notag\\
    \leq{} &CM_{k,q}(t)Q_{2,n}(t)\int_{0}^{1}\left(1-\frac{s}{t}\right)^{-\frac{n\beta}{\alpha q}}\left(\frac{s}{t}\right)^{\frac{n\beta}{\alpha q}-1-k}d\left(\frac{s}{t}\right)
\end{align}
Let $\rho=\frac{s}{t}$, we have
\begin{equation}
    t^{k}\Vert I_2\Vert_{L^{r}}\leq M_{k,q}(t)Q_{2,n}(t)\int_{0}^{1}\rho^{-\frac{n\beta}{\alpha q}}\rho^{\frac{n\beta}{\alpha q}-1-k}d\rho.
\end{equation}
Notice that $q>\frac{n}{\alpha}$ and $1k<1+\frac{n\beta}{\alpha q}$, let  $\varpi \overset{\Delta}{=}1-\frac{n\beta}{\alpha q}>0$,$\varrho\overset{\Delta}{=}\frac{n\beta}{\alpha q}-k>0$. So from the properties of the $\mathcal{B}$ function, we get
$$\int_{0}^{1}(1-\rho)^{-\frac{n\beta}{\alpha q}}\rho^{\frac{n\beta}{\alpha q}-1-k}d\rho=\mathcal{B}(\varpi,\varrho)<\infty.$$
So we have 
\begin{equation}\label{4.1.2}
    t^{k}\Vert I_2 \Vert_{L^{r}}\leq CM_{k,q}(t)Q_{2,n}(t).
\end{equation}
We next estimate $t^{k}\Vert I_3 \Vert_{L^{r}(\mathbb{R}^{n})}$, by Lemma \ref{4.1}
\begin{equation}
    \begin{split}
         t^{k}\Vert I_3 \Vert_{L^{r}(\mathbb{R}^{n})}
            &= t^{k}\left\Vert \int_{0}^{t} T_{\alpha}^{\beta}(t-s)a \tilde{u}ds \right\Vert_{L^{r}(\mathbb{R}^{n})}\\
            &\leq Ct^{k}\int_{0}^{t} (t-s)^{\beta-1} \Vert\tilde{u}\Vert_{L^{r}(\mathbb{R}^{n})}ds\\
            &\leq Ct^{k}\int_{0}^{t}t^{\beta-1}\left(1-\frac{s}{t}\right)^{\beta-1}s^{k}\Vert \tilde{u} \Vert_{L^{r}}s^{-k}ds\\
            &\leq CM_{k,q}(t)t^{\beta}\int_{0}^{1}\left(1-\frac{s}{t}\right)^{\beta-1}\left(\frac{s}{t}\right)^{-k}d\left(\frac{s}{t}\right).
    \end{split}
\end{equation}
Set $\varpi \overset{\Delta}{=}\beta>0$, $\varrho\overset{\Delta}{=}1-k>0$, and define $\rho=\frac{s}{t}$
\begin{equation}
    \int_{0}^{1}(1-\rho)^{\beta-1}\rho^{-k}d\rho=\mathcal{B}(\varpi,\varrho)<\infty.
\end{equation}
That is
\begin{equation}\label{4.1.3}
    t^{k}\Vert I_3 \Vert_{L^{r}}\leq CM_{k,q}(t)t^{\beta}.
\end{equation}

Finally estimate $t^{k}\Vert I_4 \Vert_{L^{r}(\mathbb{R}^{n})}$,
    also by Lemma \ref{4.1}, we obtain
 \begin{equation}
        \begin{split}
            t^{k}\left\Vert I_4 \right\Vert_{L^{r}(\mathbb{R}^{n})}
             &= t^{k}\Vert \int_{0}^{t} T_{\alpha}^{\beta}(t-s)b(u_{2}^{2}-u_{1}^{2})ds \Vert_{L^{r}(\mathbb{R}^{n})}\\
            &\leq Ct^{k}\int_{0}^{t} (t-s)^{\beta-1} \Vert (u_{2}^{2}-u_{1}^{2})\Vert_{L^{r}(\mathbb{R}^{n})} ds\\
            &\leq t^{k}\int_{0}^{t} (t-s)^{\beta-1} \Vert \tilde{u}(u_1+u_2)\Vert_{L^{r}(\mathbb{R}^{n})} ds\\
            &\leq Ct^{k}\int_{0}^{t} (t-s)^{\beta-1} \Vert \tilde{u}\Vert_{L^{r}(\mathbb{R}^{n})}\Vert u_1+u_2\Vert_{L^{r}(\mathbb{R}^{n})} ds\\
            &\leq Ct^{k}\int_{0}^{t} (t-s)^{\beta -1}\Vert \tilde{u}(u_1+u_2) \Vert_{L^{q}(\mathbb{R}^{n})} ds\\
            &\leq Ct^{k}\int_{0}^{t}t^{\beta-1}\left(1-\frac{s}{t}\right)^{\beta -1} s^{k}\Vert \tilde{u}\Vert_{L^{r}(\mathbb{R}^{n})} \Vert u_1+u_2 \Vert_{L^{r}(\mathbb{R}^{n})} s^{-k}ds\\
            &\leq CM_{k,q}(t) t^{\beta}\int_{0}^{1} \left(1-\frac{s}{t}\right)^{\beta -1}\left(\frac{s}{t}\right)^{-k} \Vert u_1+u_2 \Vert_{L^{r}(\mathbb{R}^{n})} d\left(\frac{s}{t}\right).
    \end{split}
    \end{equation}
We know $\Vert u_1+u_2 \Vert_{L^{r}(\mathbb{R}^{n})}\leq C.$ Let $\rho=\frac{s}{t}$, we have
\begin{equation}
    t^{k}\Vert I_4\Vert_{L^{r}(\mathbb{R}^{n})} \leq CM_{k,q}(t) t^{\beta}\int_{0}^{1}(1-\rho)^{\beta-1}\rho^{-k}d\rho.
\end{equation}
It is easy to see that $\varpi\overset{\Delta}{=}\beta>0$, $\varrho\overset{\Delta}{=}1-k>0$, then
$$\int_{0}^{1}(1-\rho)^{\beta-1}\rho^{-k}dx=\mathcal{B}(\varpi,\varrho)<\infty.$$
We have
\begin{equation}\label{4.1.4}
    t^{k}\Vert I_4 \Vert_{L^{r}(\mathbb{R}^{n})}\leq Ct^{\beta}M_{k,q}(t).
\end{equation}
Combine (\ref{4.1.1}), (\ref{4.1.2}) ,(\ref{4.1.3}) and (\ref{4.1.4}), we have
\begin{equation}
        \begin{split}
            M_{k,q}(t) 
            &= t^{k}(\Vert I_1 \Vert_{L^{r}(\mathbb{R}^{n})}+\Vert I_2 \Vert_{L^{r}(\mathbb{R}^{n})}+\Vert I_3 \Vert_{L^{r}(\mathbb{R}^{n})}+\Vert I_4 \Vert_{L^{r}(\mathbb{R}^{n})})\\
            &= C M_{k,q}(t)(Q_{1,n}(t)+Q_{2,n}(t)+2t^{\beta}).
        \end{split}
    \end{equation}
    By Lemma \ref{4.1.1.1}, we have
     $$ Q_{1,p(t)}+Q_{2,p}(t)=\underset{0<s<t}{\sup}s^{1-\frac{n\beta}{\alpha q}}\Vert u_1 \Vert_{L^{q}(\mathbb{R}^{n})} +\underset{0<s<t}{\sup}s^{1-\frac{n\beta}{\alpha q}}\Vert u_2 \Vert_{L^{q}(\mathbb{R}^{n})} \leq 2t^{\frac{\alpha\varepsilon}{n+\alpha\varepsilon}}.$$
    So there exists $t_{0}$ such that
    \begin{equation}
        M_{k,q}(t) \leq CM_{k,q}(t), {\rm for}\  \forall t\in[0,t_0)
    \end{equation}
   It can be seen from the above formula that when $t_0>0$, $M_{k,q}(t)\equiv 0$ that is, the weak solution of the system (\ref{1.1}) on $t\in [0,t_0)$ is unique. It is easy to know that if the above process is repeated with $t_0$ as the initial time, the weak solution of the system of equations on $ [t_0,2t_0)$ is unique. By repeating the above process, we get that the system of equations (\ref{1.1}) is unique on $[0,T]$. That is, the weak solution is unique in its corresponding space.\qed
\section{Blow-up criterion for weak solution}
   The main proof in this section is that if a weak solution blow up in finite time, then the all $L^{h}$-norms of the weak solution blow up at the same time for $h>q$.
\begin{thm}\label{1.1.5}
Under the same assumption as Theorem \ref{1.1.1} and $r=q+\varepsilon$ where $\varepsilon$ is small enough. Let $T_{\max}^{r}$ be the largest $L^{r}$-norm existence time of a weak solution, i.e.
\begin{align}
    \Vert u(\cdot, t)\Vert_{L^{r}(\mathbb{R}^{n})} \leq \infty,\ \ & {\rm for}\  {\rm all}\ 0<t<T_{\max}^{r}\\
    \underset{t\to T_{\max}^{r}}{\lim \sup}\Vert u(\cdot, t)\Vert_{L^{r}(\mathbb{R}^{n})}&=\infty,
\end{align}
and $T_{\max}^{r}$ be the largest $L^{h}$-norm time of a weak solution for $h\geq r>q$. Then if $T_{\max}^{r}<\infty$ for any $h$,
\begin{equation}
    T_{\max}^{h}=T_{\max}^{r},\ \ {\rm for\ } {\rm all\ } h\geq r.
\end{equation}
\end{thm}   
\begin{pro}
Since $\Vert u(\cdot,t)\Vert_{L^{1}(\mathbb{R}^{n})}\leq \Vert u_0\Vert_{L^{1}(\mathbb{R}^{n})}$, by interpolation inequality, we know that for $h\geq r$, $T_{\max}^{h}\leq T_{\max}^{r}$. If $T_{\max}^{h}<T_{\max}^{r}$ for any $h\geq r$, then we will have contradiction arguments. $T_{\max}^{h}<T_{\max}^{r}$ implies 
$$ \underset{t\to T_{\max}^{r}}{\lim \sup}\Vert u(\cdot, t)\Vert_{L^{r}(\mathbb{R}^{n})}=:A<\infty.$$
Then for $h\geq r>q$, using the interpolation inequality, Sobolev inequality and Young's inequality together, we obtain
\begin{equation}
    \begin{split}
        \Vert u \Vert_{L^{h}(\mathbb{R}^{n})}^{h}&\leq \Vert u \Vert_{L^{\frac{nh}{n-\alpha}}(\mathbb{R}^{n})}^{\frac{nh(h-r)}{nh+2r-nr}}\Vert u \Vert_{L^{r}(\mathbb{R}^{n})}^{\frac{r(\alpha n-\alpha r+2r)}{nh+2r-nr}}\\
        &= \Vert u^{\frac{h}{2}} \Vert_{L^{\frac{2n}{n-\alpha}}(\mathbb{R}^{n})}^{\frac{2n(h-r)}{nh+2r-nr}}\Vert u \Vert_{L^{r}(\mathbb{R}^{n})}^{\frac{r(\alpha n-\alpha r+2r)}{nh+2r-nr}}\\
        &\leq S_{\alpha,n}^{\frac{2n(h-r)}{nh+2r-nr}}\left\Vert D^{\frac{\alpha}{2}}u^{\frac{h}{2}} \right\Vert_{L^{2}(\mathbb{R}^{n})}^{\frac{2n(h-r)}{nh+2r-nr}}\Vert u \Vert_{L^{r}(\mathbb{R}^{n})}^{\frac{r(\alpha n-\alpha r+2r)}{nh+2r-nr}},
    \end{split}
\end{equation}
where $\frac{2n(h-r)}{nh+2r-nr}<2.$ From (\ref{3.1.9}), we can similarly obtain
\begin{equation}
    a \Vert u \Vert_{L^{h}(\mathbb{R}^{n})}^{h}\leq \hat{C_1}\left\Vert D^{\frac{\alpha}{2}}u^{\frac{h}{2}} \right\Vert_{L^{2}(\mathbb{R}^{n})}^{2}+C(h,r,a)\left(\Vert u \Vert_{L^{r}(\mathbb{R}^{n})}^{r}\right)^{\frac{\alpha(n-r)+2r}{h(n-r)+2r}}.
\end{equation}
In the same way we can get
\begin{equation}
    \begin{split}
        \Vert u \Vert_{L^{h+1}(\mathbb{R}^{n})}^{h+1}&\leq \Vert u \Vert_{L^{\frac{nh}{n-\alpha}}(\mathbb{R}^{n})}^{\frac{nh(h+1-r)}{nh+2r-nr}}\Vert u \Vert_{L^{r}(\mathbb{R}^{n})}^{\frac{r(2r-n+\alpha h-\alpha-\alpha r)}{nh+2r-nr}}\\
        &= \Vert u^{\frac{h}{2}} \Vert_{L^{\frac{2n}{n-\alpha}}(\mathbb{R}^{n})}^{\frac{2n(h+1-r)}{nh+2r-nr}}\Vert u \Vert_{L^{r}(\mathbb{R}^{n})}^{\frac{r(2r-n+\alpha h-\alpha-\alpha r)}{nh+2r-nr}}\\
        &\leq S_{\alpha,n}^{\frac{2n(h+1-r)}{nh+2r-nr}}\left\Vert D^{\frac{\alpha}{2}}u^{\frac{h}{2}} \right\Vert_{L^{2}(\mathbb{R}^{n})}^{\frac{2n(h-r)}{nh+2r-nr}}\Vert u \Vert_{L^{r}(\mathbb{R}^{n})}^{\frac{r(2r-n+\alpha h-\alpha-\alpha r)}{nh+2r-nr}},
    \end{split}
\end{equation}
where $\frac{2n(h-+1r)}{nh+2r-nr}<2.$ From (\ref{3.2.5}), we can similarly obtain
\begin{equation}
    (h-1-bh)\Vert u \Vert_{L^{h+1}(\mathbb{R}^{n})}^{h+1}\leq \hat{C_3}\left\Vert D^{\frac{\alpha}{2}}u^{\frac{h}{2}} \right\Vert_{L^{2}(\mathbb{R}^{n})}^{2}+C(h,r,b)\left(\Vert u \Vert_{L^{r}(\mathbb{R}^{n})}^{r}\right)^{\frac{2r-n+\alpha(h-1-r)}{nh+2r-nr}}.
\end{equation}
Take $\hat{C_1}=\hat{C_3}=-\frac{h-1}{h}$, then for any $h\geq r>q$, using (\ref{3.1.3}), we have
\begin{equation}
    \begin{split}
        \partial_{t}^{\alpha}\Vert u \Vert_{L^{h}(\mathbb{R}^{n})}^{h}&\leq -\frac{2(h-1)}{h}\left\Vert D^{\frac{\alpha}{2}}u^{\frac{h}{2}} \right\Vert_{L^{2}(\mathbb{R}^{n})}^{2}+C(h,r,a)\left(\Vert u \Vert_{L^{r}(\mathbb{R}^{n})}^{r}\right)^{\frac{\alpha(n-r)+2r}{h(n-r)+2r}}\\
        &+C(h,r,b)\left(\Vert u \Vert_{L^{r}(\mathbb{R}^{n})}^{r}\right)^{\frac{2r-n+\alpha(h-1-r)}{nh+2r-nr}}\\
        &\leq C(h,r,a,b,A).
    \end{split}
\end{equation}
i.e.
\begin{equation}
    \Vert u(\cdot,t)\Vert_{L^{h}(\mathbb{R}^{n})}\leq C(h,r,a,b,A, \Vert u_0\Vert_{L^{h}(\mathbb{R}^{n})}, T_{\max}^{h}), {\rm for} t\in \left(0,T_{\max}^{h}\right),
\end{equation}
which contradicts with 
$$\underset{t\to T_{\max}^{r}}{\lim \sup}\Vert u(\cdot, t)\Vert_{L^{r}(\mathbb{R}^{n})}=\infty.$$
Thus we have the conclusion that $T_{\max}^{h}=T_{\max}^{r}$ for all $h\geq r>q$, i.e. $L^{h}$-norms blow up at the same time.\qed
\end{pro}
\bibliography{Reference}
\end{document}